# STATISTICAL MECHANICAL SYSTEMS ON COMPLETE GRAPHS, INFINITE EXCHANGEABILITY, FINITE EXTENSIONS AND A DISCRETE FINITE MOMENT PROBLEM

By Thomas M. Liggett,[1] Jeffrey E. Steif[2] and Bálint Tóth[3]

*University of California at Los Angeles, Chalmers University of Technology and Göteborg University and Technical University Budapest*

We show that a large collection of statistical mechanical systems with quadratically represented Hamiltonians on the complete graph can be extended to infinite exchangeable processes. This extends a known result for the ferromagnetic Curie–Weiss Ising model and includes as well all ferromagnetic Curie–Weiss Potts and Curie–Weiss Heisenberg models. By de Finetti's theorem, this is equivalent to showing that these probability measures can be expressed as averages of product measures. We provide examples showing that "ferromagnetism" is not however in itself sufficient and also study in some detail the Curie–Weiss Ising model with an additional 3-body interaction. Finally, we study the question of how much the antiferromagnetic Curie–Weiss Ising model can be extended. In this direction, we obtain sharp asymptotic results via a solution to a new moment problem. We also obtain a "formula" for the extension which is valid in many cases.

**1. Introduction.** Let $X = (X_1, \ldots, X_n)$ be a finite exchangeable collection of random variables taking values in a space $S$ which is assumed to be a closed subset of $R^s$, often a finite set or the $(s-1)$-dimensional unit sphere. (Finite exchangeable means that the distribution is invariant under all permutations of $\{1, \ldots, n\}$.) In such a situation, one can ask whether $X$ is extendible to an infinite exchangeable process. In other words, does

Received December 2005; revised October 2006.
[1]Supported in part by NSF Grant DMS-03-01795.
[2]Supported by the Swedish Natural Science Research Council and Göran Gustafsson Foundation for Research in Natural Sciences and Medicine.
[3]Supported by the Hungarian Scientific Research Fund (OTKA) Grant T037685 and the RDSES program of the European Science Foundation.

there exist a process $Y = (Y_i)_{i \geq 1}$, taking values in $S$, whose distribution is invariant under finite permutations and such that $X$ and $(Y_1, \ldots, Y_n)$ have the same distribution? (In this case, the process $Y$ will often not be unique.) If the distribution of $X$ is "an average of product measures," meaning that the distribution of $X$ can be expressed as

$$\int_{P(S)} \mu^{\otimes n} \, d\rho(\mu),$$

where $P(S)$ is the set of probability measures on $S$, $\mu^{\otimes n}$ is the $n$-fold product of $\mu$ and $\rho$ is a probability measure on $P(S)$ (endowed with an appropriate $\sigma$-algebra), then it is immediate that $X$ is extendible to an infinite exchangeable process. One such infinite exchangeable process of course has distribution

$$\int_{P(S)} \mu^{\otimes \infty} \, d\rho(\mu). \tag{1.1}$$

The important de Finetti theorem says that any infinite exchangeable process (when $S$ is a complete separable metric space) can be expressed as in (1.1) for a unique $\rho$. Hence $X$ is extendible to an infinite exchangeable process if and only if it is "an average of product measures." When a finite exchangeable process is extendible to an infinite exchangeable process, we will write it is infinitely extendible (IE). When $S = \{0, 1\}$, $P(S)$ can be identified with $[0, 1]$ and a probability measure on $P(S)$ can be identified with a random variable $W$ taking values in $[0, 1]$. The $k$th moment, $E[W^k]$, is then the probability that the first $k$ random variables are 1. If one can extend a finite exchangeable sequence $X = (X_1, \ldots, X_n)$ to an infinite one, there may be more than one extension and so the distribution of $W$ is not unique. See [13] for some discussion concerning this point. Any $W$ which can be used will be called a *representing* $W$ for $X$.

Surveys of exchangeability can be found in [2] and [4]. The problem of determining when finite exchangeable sequences can be extended to infinite exchangeable ones has attracted some attention in the past; see the two above references as well as [21] and [23]. In the present paper, we will study this question in the mean field statistical mechanics context and also study how much one can extend if IE fails.

We start by stating a known result where $|S| = 2$. Consider the Curie–Weiss Ising model with parameters $J$ and $h$ representing the coupling constant and the external field. This is simply the Ising model on the complete graph with symmetric 1 and 2 body interactions. See [6] or [19] for an in-depth discussion of this model. This is the probability measure on $\{\pm 1\}^n$ where the probability of the configuration $\sigma$ is proportional to $e^{H(\sigma)}$ where

$$H(\sigma) = h \sum_{i=1}^{n} \sigma_i + \frac{J}{2} \left( \sum_{i=1}^{n} \sigma_i \right)^2. \tag{1.2}$$

The model is said to be *ferromagnetic* if $J \geq 0$. The following result is proved in [20] where the technique is credited to [11].

THEOREM 1.1. *For $n \geq 1$, $J \geq 0$ and any $h$, the Curie–Weiss Ising model with parameters $n$, $J$ and $h$ is IE.*

*Remarks.* When $|S| = 2$, a necessary condition for being IE is that the process be associated (i.e., increasing events are positively correlated, see Definition 2.11 on page 77 of [15]); this follows from an obvious generalization of Proposition 2.22 on page 83 of [15]. In fact, it is not hard to check that IE even implies the stronger FKG lattice condition (see (2.13) on page 78 of [15] for this definition). Using this, it is easy to check that for $n \geq 2$ and $J < 0$, the model is not IE. We point out however that an elementary example in [16] shows that this FKG lattice condition plus finite exchangeability is not sufficient for being IE. An alternative way to see that IE implies $J \geq 0$ is as follows. We have IE if and only if there is a random variable $0 \leq W \leq 1$ so that

$$(1.3) \qquad EW^k(1-W)^{n-k} = \frac{e^{h(2k-n)+(J/2)(2k-n)^2}}{\sum_{\sigma \in \{\pm 1\}^n} e^{H(\sigma)}}, \qquad 0 \leq k \leq n.$$

Hölder's inequality then gives

$$EW^k(1-W)^{n-k} \leq (EW^n)^{k/n}(E(1-W)^n)^{1-(k/n)}$$

from which it is easy to deduce that $J \geq 0$.

The following result extends Theorem 1.1 to a large number of models which have a *quadratic representation* for their Hamiltonian. It is proved using a similar method to that in [20]. Let $\nu$ be an arbitrary probability measure on $R^s$ satisfying

$$(1.4) \qquad \int_{R^s} e^{v \cdot y} \, d\nu(y) < \infty$$

for all $v \in R^s$. Next assume that $\nu$ and $n$ are such that

$$(1.5) \qquad Z_{\nu,n} := \int_{(R^s)^n} e^{(1/2) \sum_{i,j=1}^n x_i \cdot x_j} \, d\nu^{\otimes n}(x_1, \ldots, x_n) < \infty.$$

Then we can consider the probability measure $\mu_{\nu,n}$ on $(R^s)^n$ which is absolutely continuous with respect to $\nu^{\otimes n}$ with Radon–Nikodym derivative at $(x_1, \ldots, x_n)$ given by

$$\frac{e^{(1/2) \sum_{i,j=1}^n x_i \cdot x_j}}{Z_{\nu,n}}.$$

THEOREM 1.2. *For any $\nu$ and $n$ satisfying* (1.4) *and* (1.5), *$\mu_{\nu,n}$, viewed as the distribution of a finite exchangeable collection of $n$ random variables taking values in $R^s$, is IE.*

Although the above formulation is very simple, many models (including the Curie–Weiss Ising model) fall into this category as we now briefly discuss.

- Curie–Weiss Potts model:

Let $q$ be an integer larger than 1, $S = \{1, \ldots, q\}$, $J \in R$ and $h: S \to R$. The Curie–Weiss Potts model with parameters $n$, $J$ and $h$ is the probability measure on $S^n$ where the probability of the configuration $\sigma$ is proportional to $e^{H(\sigma)}$, where

$$H(\sigma) = \sum_{i=1}^n h(\sigma_i) + \frac{J}{2}\left(\sum_{i,j=1}^n I_{\{\sigma_i = \sigma_j\}}\right).$$

The model is said to be *ferromagnetic* if $J \geq 0$. (Note $q = 2$ is equivalent to the Curie–Weiss Ising model.)

- Heisenberg model:

Let $r$ be a nonnegative integer, $S$ the $r$-dimensional sphere, $J \in R$ and $h: S \to R$. Letting $dx$ denote "surface area" on $S$, the (classical) Heisenberg model with parameters $n$, $J$ and $h$ is the probability measure on $S^n$ whose Radon–Nikodym derivative with respect to $dx^{\otimes n}$ at $\sigma = (\sigma_1, \ldots, \sigma_n)$ is proportional to $e^{H(\sigma)}$, where

$$H(\sigma) = \sum_{i=1}^n h(\sigma_i) + \frac{J}{2}\left(\sum_{i,j=1}^n \sigma_i \cdot \sigma_j\right).$$

The model is said to be *ferromagnetic* if $J \geq 0$. (Note $r = 0$ is equivalent to the Curie–Weiss Ising model.)

- Curie–Weiss clock model (see [8]):

Let $q$ be an integer larger than 1, $S$ be $q$ points on the unit circle with constant spacing, $J \in R$ and $h: S \to R$. The Curie–Weiss clock model with parameters $n$, $J$ and $h$ is the probability measure on $S^n$ where the probability of the configuration $\sigma$ is proportional to $e^{H(\sigma)}$, where

$$H(\sigma) = \sum_{i=1}^n h(\sigma_i) + \frac{J}{2}\left(\sum_{i,j=1}^n \sigma_i \cdot \sigma_j\right).$$

The model is said to be *ferromagnetic* if $J \geq 0$. (Note again $q = 2$ is equivalent to the Curie–Weiss Ising model.)

COROLLARY 1.3. *For any $n$, $J \geq 0$ and $h: S \to R$, the Curie–Weiss Potts model, the Curie–Weiss Heisenberg model and the Curie–Weiss clock model are IE.*

*Remarks.* The fuzzy Potts model is obtained from the Potts model defined above by partitioning the possible spins into two sets. Häggström [10] proved that the ferromagnetic fuzzy Potts model (with no external field) on any graph has positive correlations and in fact satisfies the FKG lattice condition. We point out that for the special case of the Curie–Weiss Potts model (i.e., on the complete graph), this follows from Corollary 1.3 by using (i) the trivial fact that such a "projection" of an IE system is IE and (ii) IE systems with $|S| = 2$ satisfy the FKG lattice condition. This same argument shows that if we take the ferromagnetic Heisenberg model on the complete graph on $n$ vertices and partition the sphere into two arbitrary measurable sets, then the induced measure on $\{0,1\}^n$ satisfies the FKG lattice condition; this does not appear to be obvious directly.

It seems reasonable to ask whether all "ferromagnetic" systems on complete graphs are IE. One problem with this is that it is not clear exactly which systems should be considered ferromagnetic. We first consider the case when $|S| = 2$ but where we add 3-body interactions. Consider the probability measure $\mu_{h,J_2,J_3,n}$ on $\{\pm 1\}^n$ where the probability of the configuration $\sigma$ is proportional to $e^{H(\sigma)}$, where

$$(1.6) \qquad H(\sigma) = h \sum_{i=1}^{n} \sigma_i + \frac{J_2}{2}\left(\sum_{i=1}^{n} \sigma_i\right)^2 + \frac{J_3}{6}\left(\sum_{i=1}^{n} \sigma_i\right)^3.$$

This is the Curie–Weiss Ising model with an additional 3-body interaction term. The system is ferromagnetic (as defined in Chapter 4 of [15]) if and only if $J_2, J_3 \geq 0$. However, if $n = 3$, for example, and $h$ and $J_2$ are fixed, then for $J_3$ sufficiently large, we have

$$P(X_1 = 1 | X_2 = -1, X_3 = -1) > P(X_1 = 1 | X_2 = -1, X_3 = 1).$$

This implies that the FKG lattice condition fails and so $(X_1, X_2, X_3)$ is not IE. A similar argument works for any fixed $n \geq 3$ or for $J_3$ sufficiently negative and also shows that for any $n \geq 3$, if $J_2 = 0$ and $J_3 \neq 0$, then IE fails.

We now restrict to only 2-body interactions but general $S$. For 2-body interactions, one might define ferromagnetic to mean that the 2-body interactions are of the form $f(x \cdot y)$ where $f$ is an increasing function. However, it turns out that a system which has only 2-body interactions of this form need not be IE.

PROPOSITION 1.4. *For every $n$, there is an increasing function $f$ so that if $S = \{(1,0), (0,1), (-1,0), (0,-1)\}$ or $S = \{-1, 0, 1\}$, then the finite exchangeable probability measure on $S^n$ given by the Hamiltonian*

$$(1.7) \qquad H(\sigma) = \sum_{1 \leq i < j \leq n} f(\sigma_i \cdot \sigma_j)$$

*is not IE.*

*Remarks.* (i) The first $S$ shows that we can take the spin values to have length 1 while the second $S$ shows that we can take the spin values to be a subset of $R$.

(ii) It would be of interest to investigate whether IE would follow if one assumed that $f$ had some higher-order monotonicity.

For the rest of the results, we continue to restrict to $|S| = 2$. Unfortunately, we need to break the class of finite sequences $(X_1, \ldots, X_n)$ which are IE into two classes.

DEFINITION 1.5. If $(X_1, \ldots, X_n)$ is IE, we call it *type* 1 if there exists a random variable $W$ which represents $(X_1, \ldots, X_n)$ and is in $(0, 1)$ a.s. It is called *type* 2 otherwise [i.e., if every representing $W$ satisfies $P(W \in \{0, 1\}) > 0$].

*Remarks.* A trivial example of type 2 is where $P(X_1 \neq X_2) = 0$. In this case, a representing $W$ is trivially unique and has $P(W \in \{0, 1\}) = 1$. One can check that for $n = 2$ this is the only type 2 situation. However, for $n = 3$, we have a less trivial example of type 2 where $(X_1, X_2, X_3)$ is represented by a $W$ satisfying $P(W = 0) = P(W = 1/2) = 1/2$. In this case, if $V$ were another representing random variable (meaning that the first 3 moments are the same as those for $W$), then one can check that $E[(V^{3/2} - (1/2)V^{1/2})^2] = 0$. This implies that $V^{1/2}(V - (1/2)) = 0$ a.s. which forces $V$ to have the same distribution as $W$. It is also trivial to find type 1 $X$'s which have a representing $W$ satisfying $P(W \in \{0, 1\}) > 0$.

Concerning the problem of determining whether a given finite exchangeable process is type 1 IE, we mention the following characterization which will be used in the proofs of Propositions 1.7, 1.10 and 1.12.

PROPOSITION 1.6. *Let* $E_k = \{X_1 = \cdots = X_k = 1, X_{k+1} = 0, \ldots, X_n = 0\}$. *Then* $(X_1, \ldots, X_n)$ *is IE of type* 1 *if and only if there exists a random variable* $\xi$ *and* $c > 0$, *so that*

$$cP(E_k) = E[e^{(2k-n)\xi}]$$

*for* $k = 0, \ldots, n$.

REMARK. We observe that the latter is also a type of moment problem, since this condition is the statement that $P(E_k)$ is the $k$th moment of $e^{2\zeta}$, where $\zeta$'s distribution has Radon–Nikodym derivative $e^{-nx}$ with respect to the distribution of $\xi$.

Using Proposition 1.6, we will obtain the following result which provides some further information concerning the Curie–Weiss Ising model with an additional 3-body interaction term.

PROPOSITION 1.7. *Consider the probability measure $\mu_{h,J_2,J_3,n}$ corresponding to the Hamiltonian given in (1.6). For all $h, J_2, J_3$ with $J_3 \neq 0$, there exists $N$ such that for all $n \geq N$, $\mu_{h,J_2,J_3,n}$ is not IE.*

This should be contrasted with the fact that for all $h, J_2 > 0$ and $n$, there exists $\varepsilon > 0$ so that for $|J_3| < \varepsilon$, $\mu_{h,J_2,J_3,n}$ is IE; this follows readily from the alternative proof of Theorem 1.1 together with continuity.

Proposition 1.7 might be viewed as unnatural for the following reason. As $n$ increases, it is not so physically natural to keep the coefficients $J_2/2$ and $J_3/6$ fixed but rather they perhaps should decrease with $n$ and the appropriate Hamiltonian would be

$$(1.8) \qquad H(\sigma) = h \sum_{i=1}^n \sigma_i + \frac{J_2}{2n} \left( \sum_{i=1}^n \sigma_i \right)^2 + \frac{J_3}{6n^2} \left( \sum_{i=1}^n \sigma_i \right)^3.$$

We do not know the answer to the following question.

QUESTION 1.8. *Is it the case that for any $h$, $J_2$ and $J_3 \neq 0$, for all large $n$, the model using the Hamiltonian (1.8) is not IE?*

We have however the following two results related to this question where we take $h = 0$ for simplicity. (Proposition 1.12 below tells us that taking $h = 0$ is no restriction.)

PROPOSITION 1.9. *Consider the probability measure $\mu_{J_2,J_3,n}$ corresponding to the Hamiltonian given in (1.8) with $h = 0$. If $|J_3| > J_2$, then there exists $N$ such that for all $n \geq N$, $\mu_{J_2,J_3,n}$ is not IE.*

PROPOSITION 1.10. *Given $J_2$ and $J_3 \neq 0$, for only finitely many even values of $n$ can the system on the complete graph on $n$ vertices with Hamiltonian*

$$(1.9) \qquad H(\sigma) = \frac{J_2}{2n^2} \left( \sum_{i=1}^n \sigma_i \right)^2 + \frac{J_3}{6n^3} \left( \sum_{i=1}^n \sigma_i \right)^3$$

*be IE.*

To see in another way the degree to which the 3-body interaction term hinders being IE, we look at $n = 4$. In this case, one can check that when the Hamiltonian is taken to be

$$h \sum_{i=1}^4 \sigma_i + J_2 \sum_{1 \leq i < j \leq 4} \sigma_i \sigma_j + J_3 \sum_{1 \leq i < j < k \leq 4} \sigma_i \sigma_j \sigma_k$$

the system is IE if and only if $J_2 \geq 0$ and

$$\cosh(8J_3) \leq \cosh(4J_2) - 2e^{-8J_2}(\sinh(2J_2))^2.$$

This latter condition involving the hyperbolic functions comes from consideration of the moment characterization (2.2) of IE. From this one can conclude after some computation that if $J_2$ and $J_3$ both approach 0 with $J_2^3/J_3^2$ approaching $c$, then if $c > 1/2$, the system is eventually IE while if $c < 1/2$, the system is eventually not IE. Note the difference in the exponents (3 versus 2).

The earlier Proposition 1.6 has other applications as well.

DEFINITION 1.11. If $\mu$ is a finite exchangeable probability measure on $\{\pm 1\}^n$, we let $T_{J,h}(\mu)$ be the probability measure on $\{\pm 1\}^n$ which gives a configuration $\sigma$ probability proportional to $e^{H(\sigma)}\mu(\sigma)$, where

$$H(\sigma) = h\sum_{i=1}^n \sigma_i + \frac{J}{2}\left(\sum_{i=1}^n \sigma_i\right)^2.$$

Of course, if $\mu$ is uniform distribution, then $T_{J,h}(\mu)$ is just the Curie–Weiss Ising model. One might call $T_{J,h}(\mu)$ a "$(J,h)$-Isingization" of $\mu$.

PROPOSITION 1.12. *If the probability measure $\mu$ on $\{\pm 1\}^n$ is IE then for all $J \geq 0$ and $h \in R$, $T_{J,h}(\mu)$ is also IE.*

REMARK. This tells us that when studying the question of which models are IE, we can assume that there is no external field.

Given a finite exchangeable sequence $X = (X_1, \ldots, X_n)$ which is not IE, it is interesting to ask if it can be extended to a finite but longer exchangeable sequence.

DEFINITION 1.13. For $l > n$, a finite exchangeable sequence $X = (X_1, \ldots, X_n)$ is $l$-extendible if there exists a finite exchangeable process $Y = (Y_i)_{1 \leq i \leq l}$ such that $X$ and $(Y_1, \ldots, Y_n)$ have the same distribution. We let $E(n,l)$ denote the collection of finite exchangeable processes $X = (X_1, \ldots, X_n)$ which are $l$-extendible.

*Remarks.* It is not hard to show that $X = (X_1, \ldots, X_n)$ is IE if and only if it is $l$-extendible for every $l > n$. When $S$ is compact, this follows from an elementary compactness argument. If $S$ is not compact, a similar "tightness" argument can be easily carried out. In view of this fact, to prove the Curie–Weiss Ising model is extendible to an infinite exchangeable process, it would suffice (assuming for this discussion that $h = 0$) to show that for any $n$ and $J$ and any $l > n$, there exists $J' = J'(n, J, l)$ so that the projection of the Curie–Weiss Ising model on the complete graph of size $l$ with parameter $J'$ to $n$ vertices is the Curie–Weiss Ising model on the complete graph of size $n$ with parameter $J$. However, this is typically not true (it is however true for small $n$) and hence this approach to proving Theorem 1.1 does not work.

Proposition 1.12 says that for $J \geq 0$ and any $h$, $T_{J,h}$ leaves $E(n, \infty)$ invariant; the following is an interesting complement to this which says that this is false for finite $l > n$ even when $J = 0$.

PROPOSITION 1.14. *Given $2 \leq n < l < \infty$, there is a $\mu \in E(n, l)$ with full support and $h \in R$ such that $T_{0,h}(\mu) \notin E(n, l)$.*

REMARK. In contrast to the remark after Proposition 1.12, when we are asking about finite extensions, we cannot assume that there is no external field.

We finally consider the Curie–Weiss Ising model with $J < 0$; this is the antiferromagnetic case. *With the exception of Proposition* 1.20, *the rest of the results concern only the antiferromagnetic Curie–Weiss Ising model with parameters $J < 0$, $h$ and $n$ as defined in* (1.2). As observed following the statement of Theorem 1.1, the Curie–Weiss Ising model in not IE in this case. The following result gives us some very precise information concerning how far one can extend the model when $J$ is very close to 0. It will be convenient for our purposes to use a different parameterization for the Curie–Weiss Ising model. Given parameters $a > 0, b > 0$, we let the Gibbs state $\mu$ on $\{0,1\}^n$ have probabilities of the form $\mu\{\eta\} = a^k b^{k(n-k)}/s_n$, where $k$ is the number of 1's in the configuration $\eta$ and the normalizing constant is given by

$$s_n = \sum_{k=0}^{n} \binom{n}{k} a^k b^{k(n-k)}.$$

An easy computation shows that our two parameterizations are related by $a = e^{2h}$, $b = e^{-2J}$ and so $b \leq 1$ corresponds to $J \geq 0$. The following gives conditions under which it is $l$-extendible for large $l$. Since we will be letting $b \to 1$ and the Curie–Weiss Ising model with $b = 1$ is the product measure with density $a/(1+a)$, it is natural to denote $\rho := a/(1+a)$.

THEOREM 1.15. *Let $n \geq 2$. Let $\rho$ be as above and consider the Curie–Weiss Ising model on $\{0,1\}^n$ with $n$ and $a$ fixed, and let $b = 1 + (c/l)$.*

(a) *If $c < 1/(2\rho(1-\rho))$, then the Curie–Weiss Ising model with parameters $a$ and $b$ on $\{0,1\}^n$ is $l$-extendible for all sufficiently large $l$.*

(b) *If $c > 1/(2\rho(1-\rho))$, then the Curie–Weiss Ising model with parameters $a$ and $b$ on $\{0,1\}^n$ is not $l$-extendible for all sufficiently large $l$.*

*Remarks.* Regarding the remark after Proposition 1.14, we see here that the external field is indeed relevant to the finite extension problem and is even relevant for the $l \to \infty$ asymptotics. Note also the monotonicity in $a$ in the above result. This monotonicity however only holds in the asymptotics. Looking at the cases $n = 2, 3$ and $l = 4$ discussed at the beginning of Section 4, one sees that such monotonicity does not hold for finite $n$ and $l$. In fact, for $n = 4, l = 5$ and certain values of $a$, the set of $b$'s for which one can extend is not even an interval. We mention here that the argument involving (b) is considerably simpler than that for (a) since this part comes down to showing the nonexistence of a particular random variable by demonstrating that its first and second moments would not satisfy the Cauchy–Schwarz inequality.

*Further remarks.* One might reasonably ask how big $l$ needs to be in part (a) of Theorem 1.15. The answer undoubtedly depends in a complicated way on $n, c$ and $\rho$. For small $n$, one can say something about this using the criteria given in Section 4. If $n = 2$ and $2\rho(1-\rho)c \leq 1$ then the model is $l$-extendible for all $l \geq 2$. If $n = 3$, the answer is the same if $2\rho(1-\rho)c \leq 1$ and $\rho$ is close to $1/2$. However, if one takes $\rho \uparrow 1$ and $c \uparrow \infty$ with $2\rho(1-\rho)c = d \in (0, 1)$, then the condition for $l$-extendibility is asymptotically

$$l \geq \frac{1}{2(1-\rho)} \max\left[1, \frac{d^2}{1-d}\right].$$

One of the key steps in the analysis of the above is a solution to a new discrete finite moment problem. We first recall that determining whether an infinite sequence of numbers in $[0,1]$ is the sequence of moments of some random variable taking values in $[0,1]$ is called the Hausdorff moment problem for which a well-known sufficient and necessary condition is known. Conditions are also known which insure that a *finite* sequence of numbers in $[0,1]$ can be extended to a moment sequence. This will be used to obtain an alternative proof of Theorem 1.1. The following theorem is a solution to a discrete moment problem. It is the key to proving Theorem 1.15 and we believe it to be of independent interest.

THEOREM 1.16. *Given $v_1, \ldots, v_n$, there exists a $\{0, 1, \ldots, l\}$-valued random variable $N$ satisfying*

(1.10) $$v_k = EN^k \quad \text{for } k = 1, \ldots, n$$

*if and only if $c_0 + \sum_{i=1}^n c_i v_i \geq 0$ for all polynomials $P(x) = \sum_{i=0}^n c_i x^i$ of degree $n$ that have $n$ simple roots in $\{0, \ldots, l\}$ and are nonnegative on $\{0, \ldots, l\}$.*

*Remarks.* For the case $l = \infty$, such problems have been studied; see Chapter VII in [14]. In this sense, the above result is not such a large departure from known results. However, the technical result we use to verify the condition of Theorem 1.16 (which is contained in Theorem 4.3) is significantly different from that which appears in the treatment of earlier moment problems.

The proof of Theorem 1.15 is an existence proof and does not give a "formula" for the distribution of the extension. In Section 5, we will, for each $n$, $l$, $J < 0$ and $h$ give a formula for the distribution of the number of 1's in the extension. However, this distribution might be a signed measure in which case the formula is of course not valid. Nonetheless, if it is a distribution, then the formula will be correct. The motivation for this approach comes from trying to extend the first proof of Theorem 1.1 (either the one from [20] or the one that comes out of the proof of Theorem 1.2) to the antiferromagnetic situation.

Note that an exchangeable measure on $\{0,1\}^n$ can be identified with a probability measure on $\{0,\ldots,n\}$ which gives the distribution of the number of 1's. If $(Y_1,\ldots,Y_{M'})$ is an extension of $(X_1,\ldots,X_M)$, both exchangeable, with corresponding distributions $\mu_{M'}$ on $\{0,\ldots,M'\}$ and $\mu_M$ on $\{0,\ldots,M\}$, it is trivial to check that

$$\mu_M(m) = \sum_{m'} \frac{\binom{M}{m}\binom{M'-M}{m'-m}}{\binom{M'}{m'}} \mu_{M'}(m').$$

In this case, we say that $\mu_M$ is the *hypergeometric projection* of $\mu_{M'}$.

If $J > 0$ and $f_J$ is the density function for a normal random variable with mean 0 and variance $J$, it is easy to check that

$$Z_n := \int_{-\infty}^{\infty} (e^{(ix+h)} + e^{-(ix+h)})^n f_J(x)\,dx$$

is the normalization for the Curie–Weiss Ising model on $n$ vertices with parameters $-J$ and $h$. In the next result, for $M > 0$, $\binom{-M}{m} := (-1)^m \binom{M+m-1}{m}$.

PROPOSITION 1.17. *Fix $n$, $J > 0$, $h \geq 0$ and $l > n$. Let*

$$(1.11) \quad Q(j) := \frac{e^{h^2/(2J)}}{Z_n} \binom{l}{j} \sum_{m=0}^{\infty} \binom{-(l-n)}{m} e^{-(J/2)(2m+2l-n-2j+h/J)^2}$$

*for $j \in \{0,\ldots,l\}$. Then*

$$(1.12) \qquad \left(\sum_j \frac{\binom{n}{k}\binom{l-n}{j-k}}{\binom{l}{j}} Q(j)\right)_{0 \leq k \leq n}$$

*is the probability measure on $\{0,\ldots,n\}$ corresponding to the Curie–Weiss Ising model with parameters $n$, $-J$ and $h$. [This implies that $\sum_j Q(j) = 1$.]*

Hence if $Q(j) \geq 0$ for each $j$, then the finite extension exists and $Q$ provides a formula for the extension.

*Remarks.* We point out here that we know that there are cases, including when $h = 0$, where the extension exists but where $Q(j) < 0$ for some $j$ and hence the formula is invalid. On the other hand, if $Q(j) > 0$ for all $j$, then it is easy to see, using the fact that the hypergeometric projection is a linear mapping of full rank, that there are many other $l$-extensions besides that given by $Q$. (Proposition 1.20 is proved along these lines.)

The above theorem will help lead to the following two results and also gives us our first cases with $n < l$ where $Q \geq 0$.

PROPOSITION 1.18. *Fix $n$, $l$ and $h$. Then the Curie–Weiss Ising model with parameters $J$ and $h$ on $\{-1,1\}^n$ is $l$-extendible for all $J$ if and only if $n$ is odd, $l = n + 1$ and $h = 0$. Moreover, in this latter case, we always have $Q \geq 0$ and so provides a formula for the extension.*

REMARK. The following suggests why one might have expected part of the result stated in Proposition 1.18. If we take $h = 0$ and $J = -\infty$, then the Curie–Weiss Ising model corresponds to uniform distribution on subsets of "half" the vertices. It is clear that the projection of this distribution on even $n$ to $n - 1$ is simply this distribution on $n - 1$ while it is easy to see that this distribution on even $n$ cannot come from any distribution on $n + 1$.

PROPOSITION 1.19. *For all $n$, there exists $\varepsilon > 0$ such that for all $J$ with $|J| < \varepsilon$, the Curie–Weiss Ising model with parameters $h = 0$ and $J$ on $\{-1,1\}^n$ is $n + 1$-extendible and $Q \geq 0$ and so provides a formula for the extension.*

REMARK. The key part of the above result is that this is another case in which (1.11) yields a formula for the extension. The weaker fact that the system is $n + 1$-extendible for small $J$ follows from a much more general result. This is stated in the following easy proposition.

PROPOSITION 1.20. *Fix $n < l < \infty$ and let $X = (X_1, \ldots, X_n)$ be i.i.d. $0, 1$ valued with $P(X_1 = 1) = p \in (0, 1)$. Then any small perturbation of $X$ which is in $E(n, n)$ also lies in $E(n, l)$. This is false if $p = 0$ or $1$.*

Our final theorem gives us a very large number of cases where our formula is valid. These cases are qualitatively similar to those covered in Theorem 1.15 but quantitatively different.

THEOREM 1.21. *Let $n \geq 1$ and $c > 0$. Define*

$$\alpha(c) := c + \cosh^{-1}\left(\frac{1}{(1-e^{-1/c})^{1/2}}\right) - c\tanh\left[\cosh^{-1}\left(\frac{1}{(1-e^{-1/c})^{1/2}}\right)\right],$$

$$\beta(c) := \ln(\sqrt{c} + \sqrt{c-1}) + c - \sqrt{c^2 - c} \qquad (\text{defined only for } c > 1)$$

*and*

$$h^*(c) := \begin{cases} \max\{c, \alpha(c)\}, & \text{if } c \in (0,1], \\ \max\{\alpha(c), \beta(c)\}, & \text{if } c \in (1, 3/2), \\ \beta(c), & \text{if } c \geq 3/2. \end{cases}$$

*Then, for all $c > 0$ and $h > h^*(c)$, if $J = c/l$, $l$ is sufficiently large and $Q$ is defined as in* (1.11), *we have that $Q \geq 0$ and so provides a formula for the relevant extension.* [*Recall the hypergeometric projection of $Q$ is the Curie–Weiss Ising model with parameters $n, -J$ and $h$ as defined in* (1.2).]

*Remarks.* The parameterizations here and in Theorem 1.15 are different and it is easy to check that a value of $c$ in Theorem 1.15 corresponds to $2c$ in the above result. In Theorem 1.15, it is easily checked that $c < 1/(2\rho(1-\rho))$ provided that (i) $c < 2$ (no matter what $h$ is) or (ii) $c \geq 2$ and $|h| > \ln(\sqrt{c/2} + \sqrt{c/2 - 1})$. Hence, in view of the above relationship between the $c$'s in the two results and the fact that $\lim_{c\to\infty} |\beta(c) - \ln(\sqrt{c} + \sqrt{c-1})| = 1/2$, we have that for large $c$ (equivalently large $h$), the two conditions are not too far apart. We are not sure whether to believe that in the asymptotic regime where Theorem 1.15 guarantees an extension we also have that $Q \geq 0$. At the same time, we are quite sure that the bounds given in this theorem are not sharp. In this theorem, it is also possible to take $n$ growing to $\infty$ as long as $n = o(\sqrt{l})$ but we do not bother to elaborate on this.

The rest of the paper is organized as follows. In Section 2, we first give proofs of Theorem 1.2 as well as Corollary 1.3. We then give an alternative proof of Theorem 1.1. For this alternative proof, we will use a known sufficient condition for when a given finite sequence extends to an infinite moment sequence. Often such conditions, while being of theoretic interest, are difficult to apply; we find it interesting that such a condition can be checked in this concrete situation. We point out that the alternative proof of Theorem 1.1 is not as unrelated to the proof of Theorem 1.2 (in the case of the Ising model) as may at first seem. In the first proof, one explicitly uses the fact that $e^{x^2/2}$ is a moment generating function. Had we not known this, we would have to check it by verifying certain conditions, which would probably involve positivity of certain determinants which is the approach of the alternative proof. In addition, this alternative proof is the basis of the proof of Theorem 1.15. Finally, we give the proof of Proposition 1.4 in this section. In Section 3, we will give the proofs of Propositions 1.6, 1.7, 1.9, 1.10, 1.12 and 1.14. In Section 4, the proofs of Theorems 1.15 and 1.16 will be given. Finally, in Section 5, we prove Propositions 1.17, 1.18, 1.19 and 1.20 as well as Theorem 1.21.

**2. Hamiltonians with quadratic representation are IE.** We first prove Theorem 1.2.

PROOF OF THEOREM 1.2. Let $V$ be a random variable in $R^s$ whose density (with respect to $s$-dimensional Lebesgue measure) is given by

$$\rho(w) = \frac{(\int_{R^s} e^{w \cdot x} \, d\nu(x))^n}{Z_{\nu,n}} \frac{e^{-\|w\|^2/2}}{(2\pi)^{s/2}},$$

where $Z_{\nu,n}$ is given in (1.5). (It is easy to check that $Z_{\nu,n}$ is the correct normalization to yield a probability density.) For $w \in R^s$, let $P_w$ be the probability measure on $R^s$ which is absolutely continuous with respect to $\nu$ with Radon–Nikodym derivative

$$\frac{dP_w(x)}{d\nu(x)} := \frac{e^{w \cdot x}}{\int_{R^s} e^{w \cdot y} \, d\nu(y)}.$$

[Note the denominator is finite by (1.4).] Letting $F_V$ be the distribution of $V$, we claim that

$$\mu_{\nu,n} = \int_{R^s} P_w^{\otimes n} \, dF_V(w)$$

which proves the result. To see this, note that the right-hand side clearly is absolutely continuous with respect to $\nu^{\otimes n}$ with Radon–Nikodym derivative at $(x_1, \ldots, x_n)$ given by

$$\int_{R^s} \frac{\prod_{i=1}^n e^{w \cdot x_i}}{(\int_{R^s} e^{w \cdot y} \, d\nu(y))^n} \, dF_V(w) = \frac{1}{Z_{\nu,n}(2\pi)^{s/2}} \int_{R^s} e^{w \cdot \sum_{i=1}^n x_i} e^{-\|w\|^2/2} \, dw$$

which, using the formula for the moment generating function for a multidimensional standard normal random variable, is

$$\frac{1}{Z_{\nu,n}} e^{(1/2)\|\sum_{i=1}^n x_i\|^2} = \frac{1}{Z_{\nu,n}} e^{(1/2) \sum_{i,j=1}^n x_i \cdot x_j}$$

as desired. □

PROOF OF COROLLARY 1.3. We assume throughout $J > 0$ as otherwise the results are trivial.

Curie–Weiss Potts model: Let the parameters $q, J > 0$ and $h$ be given. Let $s = q$, choose $q$ orthogonal vectors $Q = \{a_1, \ldots, a_q\}$ on the sphere $\{x : \|x\| = \sqrt{J}\}$ in $R^q$ and let $\nu$ be the probability measure on $R^q$ concentrated on $Q$ giving $a_i$ weight proportional to $e^{h(i)}$. It is easy to check that the corresponding measure $\mu_{\nu,n}$ is exactly the Curie–Weiss Potts model where $a_i$ is identified with $i$.

Heisenberg model: Let the parameters $r$, $J > 0$ and $h$ be given. Let $s = r + 1$, consider the sphere $\{x : \|x\| = \sqrt{J}\}$ in $R^s$ and let $\nu$ be the probability measure on $R^s$ concentrated on this sphere whose density with respect to surface area is proportional to $e^{h(x/\sqrt{J})}$. It is easy to check that the corresponding measure $\mu_{\nu,n}$ is exactly the Heisenberg model with the sphere $\{x : \|x\| = \sqrt{J}\}$ trivially identified with the unit sphere.

Curie–Weiss clock model: Let the parameters $q$, $J > 0$ and $h$ be given. Let $s = 2$ and choose $q$ equally spaced points $Q = \{a_1, \ldots, a_q\}$ on the circle $\{x : \|x\| = \sqrt{J}\}$. Let $\nu$ be the probability measure on $R^2$ concentrated on $Q$ giving $a_i$ weight proportional to $e^{h(i)}$. It is easy to check that the corresponding measure $\mu_{\nu,n}$ is exactly the Curie–Weiss clock model where $a_i$ is identified with $i$. $\square$

We now give an alternative proof of Theorem 1.1. For this proof, we prefer to use the second parameterization given right after Proposition 1.14.

ALTERNATIVE PROOF OF THEOREM 1.1. Define sequences $u_k$ and $v_k$ for $0 \leq k \leq n$ by

$$u_k = a^k b^{k(n-k)}/s_n$$

and

$$v_k = \sum_{j=0}^{n-k} \binom{n-k}{j} u_{k+j}.$$

Then (1.3) is equivalent to

(2.1) $$EW^k = v_k, \qquad 0 \leq k \leq n.$$

Thus our problem is reduced to determining whether $v_0, \ldots, v_n$ can be extended to the sequence of moments of a random variable $0 \leq W \leq 1$. There is a classical solution to this problem. See [22] or [24] for a more recent discussion. The condition for extendibility is a bit different depending on the parity of $n$, so we assume from now on that $n$ is even, and write $n = 2m$. The odd case is similar.

A sufficient condition for the existence of $0 \leq W \leq 1$ satisfying (2.1) is that the following two matrices be (strictly) positive definite:

(2.2) $\begin{pmatrix} v_0 & v_1 & \cdots & v_m \\ v_1 & v_2 & \cdots & v_{m+1} \\ \cdots & \cdots & \cdots & \cdots \\ v_m & v_{m+1} & \cdots & v_n \end{pmatrix}$ and $\begin{pmatrix} w_1 & w_2 & \cdots & w_m \\ w_2 & w_3 & \cdots & w_{m+1} \\ \cdots & \cdots & \cdots & \cdots \\ w_m & w_{m+1} & \cdots & w_{n-1} \end{pmatrix},$

where $w_k = v_k - v_{k+1}$. We will consider only the first of these, since the treatment of the second is similar.

A matrix is positive definite if and only if its principal minors are all strictly positive. This is a standard result in linear algebra. It is usually proved by induction. However, in Proposition 4.2, we will have a result that yields an immediate proof of this fact. Usually, one thinks of the principal minors as being the determinants of the submatrices that are situated in the upper left corner of the matrix. However, it is much more convenient for our purpose to consider the ones that are situated in the lower right corner of the matrix. We will use the following notation:

$$f(k,l) = \sum_{j=0}^{k} \binom{k}{j} u_{l-j}, \qquad 0 \leq k \leq l \leq n.$$

Thus $v_{n-k} = f(k,n)$. With this notation, we must prove the strict positivity of

$$(2.3) \quad \begin{vmatrix} v_{n-2k} & \cdots & v_{n-k} \\ \cdots & \cdots & \cdots \\ v_{n-k} & \cdots & v_n \end{vmatrix} = \begin{vmatrix} f(2k,n) & \cdots & f(k,n) \\ \cdots & \cdots & \cdots \\ f(k,n) & \cdots & f(0,n) \end{vmatrix}, \qquad k = 0, \ldots, m.$$

We begin by performing some row operations. We will use repeatedly the easily verified relation $f(k,l) - f(k-1,l) = f(k-1,l-1)$. Now, subtract the second row from the first, then the third from the second, ..., and finally the last row from the $k$th row. The result is that (2.3) equals

$$\begin{vmatrix} f(2k-1,n-1) & f(2k-2,n-1) & \cdots & f(k-1,n-1) \\ \cdots & \cdots & \cdots & \\ f(k,n-1) & f(k-1,n-1) & \cdots & f(0,n-1) \\ f(k,n) & f(k-1,n) & \cdots & f(0,n) \end{vmatrix}.$$

Repeat this process a number of times, each time using one less row than the previous time. The result is that (2.3) equals

$$\begin{vmatrix} f(k,n-k) & \cdots & f(0,n-k) \\ \cdots & \cdots & \cdots \\ f(k,n) & \cdots & f(0,n) \end{vmatrix}.$$

Finally, repeat this whole procedure using columns instead of rows. The result is that (2.3) equals

$$(2.4) \quad \begin{vmatrix} f(0,n-2k) & \cdots & f(0,n-k) \\ \cdots & \cdots & \cdots \\ f(0,n-k) & \cdots & f(0,n) \end{vmatrix} = \begin{vmatrix} u_{n-2k} & \cdots & u_{n-k} \\ \cdots & \cdots & \cdots \\ u_{n-k} & \cdots & u_n \end{vmatrix}.$$

Note that up to this point, we have not used the particular form of the $u_k$'s. The statement that the left-hand side of (2.3) equals the right-hand side of (2.4) holds for any finite exchangeable measure.

Now we do use the particular form for the $u_k$'s, and assume that $b < 1$, since when $b = 1$, $\eta_1, \ldots, \eta_n$ are i.i.d., so its extendibility is immediate. Noting that one can factor out powers of $s_n$, $a$ and $b$, write the right-hand side of (2.4) as

$$\frac{1}{s_n^{k+1}} \begin{vmatrix} a^{n-2k}b^{2k(n-2k)} & \cdots & a^{n-k}b^{k(n-k)} \\ \cdots & \cdots & \cdots \\ a^{n-k}b^{k(n-k)} & \cdots & a^n b^0 \end{vmatrix}$$

$$= \frac{a^{(k+1)(n-k)}}{s_n^{k+1}} \begin{vmatrix} b^{2k(n-2k)} & \cdots & b^{k(n-k)} \\ \cdots & \cdots & \cdots \\ b^{k(n-k)} & \cdots & b^0 \end{vmatrix}$$

$$= \frac{a^{(k+1)(n-k)} b^{k(k+1)(3n-2k-1)/3}}{s_n^{k+1}} \begin{vmatrix} b^{-2k^2} & \cdots & b^{-4k} & b^{-2k} & 1 \\ \cdots & \cdots & \cdots & \cdots & \cdots \\ b^{-4k} & \cdots & b^{-8} & b^{-4} & 1 \\ b^{-2k} & \cdots & b^{-4} & b^{-2} & 1 \\ 1 & \cdots & 1 & 1 & 1 \end{vmatrix}.$$

This last determinant is of Vandermonde form, so can be computed explicitly as

$$\prod_{0 \leq j < l \leq k} (b^{-2l} - b^{-2j}),$$

which is strictly positive if $b < 1$. □

We lastly give the proof of Proposition 1.4.

PROOF OF PROPOSITION 1.4. Fix $n \geq 2$. Let $f$ satisfy $f(1) = a_n$, $f(0) = 0$ and $f(-1) = -b_n$ where $a_n$ is positive and small and where $b_n$ is positive and large.

Case 1: $S = \{(1,0), (0,1), (-1,0), (0,-1)\}$.

If the model were IE, then so would be the process on $\{0,1\}^n$, called $(Y_1, \ldots, Y_n)$, obtained by partitioning $S$ into $\{(1,0), (-1,0)\}$ and $\{(0,1), (0,-1)\}$ and letting the first set correspond to 1 and the second set to 0. This latter measure would then satisfy the FKG lattice condition which we now show it does not. It is easy to check that if $a_n$ is sufficiently small and $b_n$ is sufficiently large, then for any $(u_2, \ldots, u_n) \in \{(1,0), (-1,0)\}^{n-1}$,

$$P(Y_1 = 1 | X_2 = u_2, \ldots, X_n = u_n) < 1/2$$

which implies

$$P(Y_1 = 1 | Y_2 = 1, \ldots, Y_n = 1) < 1/2.$$

By symmetry we have

$$P(Y_1 = 1 | Y_2 = 0, \ldots, Y_n = 0) > 1/2,$$

which then violates the FKG lattice condition.

Case 2: $S = \{-1, 0, 1\}$.

Consider the process on $\{0,1\}^n$, called $(Y_1, \ldots, Y_n)$, obtained by partitioning $S$ into $\{1, -1\}$ and $\{0\}$ and letting the first set correspond to 1 and the second set to 0. It is clear that

$$P(Y_1 = 1 | Y_2 = 0, \ldots, Y_n = 0) = 2/3$$

for any choice of $a_n$ and $b_n$. On the other hand, it is easy to check that if $a_n$ is sufficiently small and $b_n$ is sufficiently large, then

$$P(Y_1 = 1 | Y_2 = 1, \ldots, Y_n = 1) < 0.51$$

leading to a similar contradiction as in case 1. □

### 3. Moment generating functions, 3-body interactions and Isingization.

3.1. *Moment generating functions and type* 1 *IE.* In this subsection, we prove Proposition 1.6.

PROOF OF PROPOSITION 1.6. For the "if" direction, let $\xi$ be as given and let $Y$ be the random variable whose distribution is absolutely continuous with respect to the distribution of $\xi$ with Radon–Nikodym derivative given by

$$\frac{(2 \cosh x)^n}{b},$$

where $b$ is $E[(2 \cosh \xi)^n]$. Note $b < \infty$ since $E[e^{(2k-n)\xi}] < \infty$ for $k = 0, \ldots, n$ and so $Y$ is well defined. Let $W = e^Y / (2 \cosh Y)$ [which is in $(0,1)$]. Then for $k = 0, \ldots, n$,

$$E[W^k(1-W)^{n-k}] = E\left[\frac{e^{kY} e^{-(n-k)Y}}{(2\cosh Y)^n}\right] = \frac{E[e^{k\xi} e^{-(n-k)\xi}]}{b} = \frac{cP(E_k)}{b}.$$

[Clearly the last term is then just $P(E_k)$, concluding that $c = b$.] This shows $\mu$ is extendible to an infinite exchangeable process with a mixing variable $W$ a.s. contained in $(0,1)$.

The only if direction is more or less obtained by going backward. Choose a random variable $W$ contained in $(0,1)$ a.s. such that

$$E[W^k(1-W)^{n-k}] = P(E_k)$$

for $k = 0, \ldots, n$. Let $Y$ be the random variable defined by $W = \frac{e^Y}{2\cosh Y}$ [here $W \in (0,1)$ is being used]. Let $\xi$ be the random variable whose distribution is absolutely continuous with respect to the distribution of $Y$ with Radon–Nikodym derivative given by

$$\frac{1}{b(2\cosh x)^n},$$

where $b$ is $E[\frac{1}{(2\cosh Y)^n}]$. Since cosh is bounded away from 0, $\xi$ is well defined. We then have for $k = 0, \ldots, n$

$$E[e^{(2k-n)\xi}] = E\left[\frac{e^{(2k-n)Y}}{b(2\cosh Y)^n}\right] = \frac{E[W^k(1-W)^{n-k}]}{b} = \frac{P(E_k)}{b}.$$

(In this case, $b = 1/c$.) □

REMARK. It can be shown that in Theorem 1.1, the $\xi$ in the above result is simply a normal random variable with mean $h$ and variance $J$.

3.2. *Curie–Weiss Ising model with 3-body interactions.* In this subsection, we prove Propositions 1.7, 1.9 and 1.10.

PROOF OF PROPOSITION 1.7. Fix $h$, $J_2$ and $J_3 \neq 0$. Choose $N$ so that for all $n \geq N$, the function

$$f := hx + \frac{J_2}{2}x^2 + \frac{J_3}{6}x^3$$

defined at the points $x = -n, -n+2, \ldots, n-2, n$ does not extend to a convex function on $[-n, n]$. Such an $N$ clearly exists by looking at $-n, -n+2, -n+4$ if $J_3 > 0$ and at $n-4, n-2, n$ if $J_3 < 0$. Fix $n \geq N$. We claim that $\mu_{h,J_2,J_3,n}$ is not IE. Choose $\varepsilon > 0$ so that any function $g$ defined at the points $x = -n, -n+2, \ldots, n-2, n$ satisfying $|f(x) - g(x)| < \varepsilon$ for each such $x$ does not have a convex extension to $[-n, n]$. If $\mu_{h,J_2,J_3,n}$ is IE, let $W$ be a representing random variable. By perturbing $W$ a little bit, we can obtain a random variable $W'$ taking values in $(0,1)$ with

(3.1) $\quad |\log[\mu(E_k)] - \log[\mu_{h,J_2,J_3,n}(E_k)]| < \varepsilon \quad$ for $k = 0, \ldots, n$,

where $\mu$ is the probability measure on $\{\pm 1\}^n$ coming from the mixing variable $W'$ and $E_k$ is as in Proposition 1.6. (This uses the fact that $\mu_{h,J_2,J_3,n}$ has full support.)

Since $W'$ takes values in $(0,1)$, Proposition 1.6 tells us that there is a random variable $\xi$ and $c > 0$ satisfying

$$E[e^{(2k-n)\xi}] = c\mu(E_k)$$

for $k = 0, \ldots, n$. Since a moment generating function exists on an interval and its logarithm is convex, we conclude that the function $h(t) := \log(E[e^{t\xi}])$ exists and is convex on $[-n, n]$. Note that its restriction to $x \in \{-n, -n+$ $2,\ldots,n-2,n\}$ differs from the function $\log[\mu(E_{(x+n)/2})]$ by a constant. In view of (3.1) and the definition of $\mu_{h,J_2,J_3,n}$, we conclude that $|h(x) - f(x)| < \varepsilon$ (after a translation of $f$ or $h$) for $x \in \{-n, -n+2, \ldots, n-2, n\}$. This is a contradiction. □

PROOF OF PROPOSITION 1.9. Since $J_3 > J_2$, it is easy to check that for large $n$ the function

$$f := \frac{J_2}{2n}x^2 + \frac{J_3}{6n^2}x^3$$

defined at the points $x = -n, -n+2, \ldots, n-2, n$ does not extend to a convex function on $[-n, n]$. From here, one can simply carry out the proof of Proposition 1.7. □

PROOF OF PROPOSITION 1.10. Fix $J_2$ and $J_3 \neq 0$ and denote the relevant measure on $\{\pm 1\}^n$ by $\mu_n$ (ignoring explicit notation of the dependence on $J_2$ and $J_3$). For simplicity, we can assume that for all even $n$, $\mu_n$ is IE. Fix $n$. If $\mu_n$ were of type 1, there would exist, by Proposition 1.6, a random variable $X_n$ such that

$$(3.2)\quad E[e^{kX_n}] = e^{(J_2/(2n^2))k^2 + (J_3/(6n^3))k^3} \qquad \text{for } k = -n, -n+2, \ldots, n-2, n.$$

Note Proposition 1.6 only says that the left- and right-hand sides are proportional but by taking $k = 0$, we see that equality holds (this is why we take $n$ even). Since $\mu_n$ need not be of type 1, we need to make a preliminary detour to (almost) obtain (3.2). We first find a type 1 IE measure $\mu_{n,m}$ on $\{\pm 1\}^n$ with $\|\mu_{n,m} - \mu_n\| < 1/m$ where $\|\cdot\|$ is total variation norm. [This can be easily done by taking a representing random variable $W$ for $\mu_n$ and perturbing it a small bit obtaining a random variable $W'$ taking values in $(0, 1)$.] By Proposition 1.6, there is a random variable $X_{n,m}$ and $c_{n,m} > 0$ such that

$$(3.3)\quad E[e^{kX_{n,m}}] = c_{n,m}\mu_{n,m}(E_{(k+n)/2}) \qquad \text{for } k = -n, -n+2, \ldots, n-2, n,$$

where $E_k$ is as in Proposition 1.6. Setting $k = 0$, we see that $\lim_{m\to\infty} c_{n,m} = [\mu_n(E_{n/2})]^{-1}$. It follows that

$$\lim_{m\to\infty} E[e^{kX_{n,m}}] = \frac{\mu_n(E_{(k+n)/2})}{\mu_n(E_{n/2})} = e^{(J_2/(2n^2))k^2 + (J_3/(6n^3))k^3}$$

(3.4)

$$\text{for } k = -n, -n+2, \ldots, n-2, n.$$

Since

$$(3.5)\quad \lim_{m\to\infty} E[e^{nX_{n,m}}] = e^{J_2/2 + J_3/6}, \qquad \lim_{m\to\infty} E[e^{-nX_{n,m}}] = e^{J_2/2 - J_3/6},$$

it follows that $\{X_{n,m}\}_{m\geq 1}$ is tight. We conclude that for some $m_\ell \to \infty$, $X_{n,m_\ell} \to X_n$ in distribution and (3.5) allows us to conclude (by uniform integrability) that

$$(3.6) \quad \lim_{m\to\infty} E[e^{kX_{n,m_\ell}}] = E[e^{kX_n}] \qquad \text{for } k = -n+2,\ldots,n-2.$$

Now, (3.4) finally allows us to conclude that

$$(3.7) \quad E[e^{kX_n}] = e^{(J_2/(2n^2))k^2 + (J_3/(6n^3))k^3} \qquad \text{for } k = -n+2,\ldots,n-2$$

which is only slightly weaker than (3.2).

Equation (3.7) now tells us that the sequence $\{nX_n\}_{n\geq 1}$ is tight and hence converges along a subsequence to $X_\infty$. For $z \in (-1,1)$, one can let $k_n/n$ approach $z$ (with $|k_n| \leq n-2$) and conclude (using $|z| < 1$ implies uniform integrability) that

$$E[e^{zX_\infty}] = e^{(J_2/2)z^2 + (J_3/6)z^3}, \qquad z \in (-1,1).$$

The two sides are complex analytic functions in $\{z : |\mathrm{Re}(z)| < 1\}$ and hence agree on the imaginary axis. It follows that

$$E[e^{itX_\infty}] = e^{(-J_2/2)t^2 - i(J_3/6)t^3}$$

for all $t \in R$. We claim that there is no random variable with this characteristic function when $J_3 \neq 0$. If there were, let $X_1$ and $X_2$ be independent copies with this distribution and we conclude that $X_1 - X_2$ is normal. Theorem 19 of [3] however says that if a sum of two independent random variables is normally distributed, then so is each summand. This yields a contradiction. □

### 3.3. Results for Isingization.
In this final subsection, we prove Propositions 1.12 and 1.14.

PROOF OF PROPOSITION 1.12. We break the proof up into two steps. We first prove the result under the further assumption that $\mu$ is of type 1 using Proposition 1.6 and then extend the result in general.

Fix a probability measure $\mu$ on $\{\pm 1\}^n$ which is type 1 IE and let $J \geq 0$ and $h \in R$. Let $E_k$ be the event that exactly the first $k$ variables are 1. By Proposition 1.6, there exists a random variable $\xi$ and $c$ so that

$$E[e^{(2k-n)\xi}] = c\mu(E_k)$$

for $k = 0, \ldots, n$. Since $T_{J,h}(\mu)(E_k)$ is proportional to $e^{h(2k-n)+(J/2)(2k-n)^2}\mu(E_k)$, it follows that $[T_{J,h}(\mu)](E_k)$ is proportional to $e^{h(2k-n)+(J/2)(2k-n)^2}E[e^{(2k-n)\xi}]$. However, the latter is clearly equal to $E[e^{(2k-n)(\xi+h+\sqrt{J}U)}]$ where $U$ is a standard normal random variable independent of $\xi$. By Proposition 1.6, we conclude that $T_{J,h}(\mu)$ is IE (and of type 1 although we do not need that).

For the second step, consider a probability measure $\mu$ on $\{\pm 1\}^n$ which is IE (perhaps of type 2). Let $W$ be a representing mixing random variable. If $P(W \in \{0,1\}) = 1$, the result is trivial. Otherwise let $m$ be the conditional distribution of $W$ given $W \notin \{0,1\}$ and let $\nu$ be the probability measure on $\{\pm 1\}^n$ given by a mixing random variable having distribution $m$. After some reflection, one sees that $T_{J,h}(\mu)$ is a convex combination of $T_{J,h}(\nu)$, $\delta_1$ and $\delta_{-1}$ where $\delta_i$ is the measure concentrating on having only $i$'s. By the first part, we know that $T_{J,h}(\nu)$ is IE and so we can conclude that $T_{J,h}(\mu)$ is IE as well. □

PROOF OF PROPOSITION 1.14. Fix $2 \leq n < l$. $E(l,l)$ can be identified with probability vectors $(g_0, \ldots, g_l)$ of length $l+1$ where $g_i$ is the probability of having $i$ 1's. Define such a probability vector by letting $g_0 = \delta$, $g_i = 0$ for $i = 1, \ldots, l-n$ and $g_i = (1-\delta)/n$ for $i = l-n+1, \ldots, l$ where $\delta > 0$ will be determined later. Let now $\mu$ be the distribution of the first $n$ variables of the $l$ finite exchangeable random variables corresponding to this $g$. It is easy to check that $\mu$ has full support. We now claim that for $h$ sufficiently close to $-\infty$ and $\delta$ sufficiently small, $T_{0,h}(\mu) \notin E(n,l)$.

Let $O := \{|\{i : X_i = 1\}| = 1\}$. Let $\nu_i$ be the probability measure in $E(n,l)$ which comes from the probability measure in $E(l,l)$ corresponding to the $g$ with $g_i = 1$. Since $2 \leq n < l$, it is clear that $\nu_i(O) < 1$ for each $i$ and so we can choose $\varepsilon > 0$ so that $\nu_i(O) \leq 1 - \varepsilon$ for each $i = 0, \ldots, l$. Since the natural map from $E(l,l)$ to $E(n,l)$ is affine, we conclude that $\nu(O) \leq 1 - \varepsilon$ for any $\nu \in E(n,l)$. However, it is clear we can choose $h$ sufficiently negative and $\delta$ sufficiently small so that $T_{0,h}(\mu)(O) > 1 - \varepsilon$ implying $T_{0,h}(\mu) \notin E(n,l)$. □

*Remarks.* $\delta$ is only used to yield full support of $\mu$. In the above proof, the $\nu$ in $E(l,l)$ whose projection was the desired $\mu$ did not have full support; by a small perturbation, one can also take this $\nu$ to have full support. There is another explanation of why the above result is true. $E(n,n)$ is foliated by differentiable curves of the form $(T_{0,h}(\mu))_{h \in R}$ as $\mu$ varies over $E(n,n)$ and at the same time, $E(n,n)$ is an $n$-dimensional simplex. This however is not a contradiction since the corner points of this simplex are fixed by $T_{0,h}$ and so these points reduce to curves having only 1 point. However, the polytope $E(n,l)$ has many corner points which are not fixed by $T_{0,h}$ and the same is true nearby these points. If $T_{0,h}$ left $E(n,l)$ invariant, then $E(n,l)$ would be foliated by regular differentiable curves near these corner points which of course cannot happen.

**4. A discrete moment problem and finite extendibility.** In this section, we prove Theorems 1.15 and 1.16.

Finding necessary and sufficient conditions for $l$-extendibility for general $n$ and $l$ seems to be out of the question, since they would necessarily be very complex. We illustrate this by stating the necessary and sufficient conditions for small $n$ and $l$. Checking these involves routine but somewhat tedious computations. One simply writes down the distribution of a general exchangeable measure on $\{0,1\}^l$, and solves the equations that guarantee that the $n$-dimensional marginals are the ones that correspond to the Curie–Weiss Ising model. Then one determines conditions guaranteeing the feasibility of the resulting linear programming problem. Of course by symmetry, one can go between the cases $a < 1$ and $a > 1$ by replacing $a$ by $1/a$ and so we just take $a \leq 1$ throughout. Here are the results:

$n = 2, l = 3$:
$$b \leq a + a^{-1} \qquad \text{if } a \leq 1,$$

$n = 2, l = 4$:
$$b \leq \tfrac{3}{2}a + \tfrac{1}{2}a^{-1} \qquad \text{if } a \leq 1,$$

$n = 2, l = 5$:
$$b \leq \begin{cases} 2a + \tfrac{1}{3}a^{-1}, & \text{if } a \leq 1/\sqrt{3}, \\ \tfrac{3}{4}a + \tfrac{3}{4}a^{-1}, & \text{if } 1/\sqrt{3} \leq a \leq 1, \end{cases}$$

$n = 3, l = 4$:
$$b \leq 1/\sqrt{a(1-a)} \qquad \text{if } a < 1,$$

$n = 3, l = 5$:
$$b \leq \begin{cases} 1/\sqrt{a(2-3a)}, & \text{if } a \leq 1/2, \\ \sqrt{3}/\sqrt{a(2-a)}, & \text{if } 1/2 \leq a \leq 1, \end{cases}$$

$n = 4, l = 5$:
$$b \in \begin{cases} (0, b_1(a)] \cup [b_2(a), a + a^{-1}], & \text{if } a \leq a_0, \\ (0, a + a^{-1}], & \text{if } a_0 \leq a \leq 1, \end{cases}$$

where $a_0 = 0.477\ldots$ is a root of $27a^8 - 148a^6 + 162a^4 - 148a^2 + 27 = 0$ and $b_1(a) \leq b_2(a)$ are the two real roots of $a^4 - a^3b^3 + a^2b^4 - ab^3 + 1 = 0$. Note that in the last case, the set of $b$'s for which one can extend is not even an interval. Note also the difference between $n = 3, l = 4$ and $n = 4, l = 5$ when $a = 1$; in the first case $b$ can be arbitrarily large, while in the second, it cannot be. These are of course just special cases of Proposition 1.18.

For a more systematic approach, one can try to imitate the alternative proof of Theorem 1.1. That proof involved two elements: (a) de Finetti's theorem, which reduces the extendibility problem to a moment problem, and (b) the solution to the moment problem given by (2.2).

The analogue of de Finetti's theorem for finite extendibility is elementary and well known. (See page 536 of [7], e.g.) Suppose that $Y_1, \ldots, Y_l$ are exchangeable $\{0,1\}$-valued random variables and set $N = \sum_{i=1}^{l} Y_i$. Then for $k \leq l$

$$P(Y_1 = 1, \ldots, Y_k = 1) = \sum_{i=k}^{l} \binom{l-k}{i-k} \frac{P(N=i)}{\binom{l}{i}} = \frac{EN(N-1)\cdots(N-k+1)}{l(l-1)\cdots(l-k+1)}.$$

Thus the $\{0,1\}$-valued exchangeable sequence $X_1, \ldots, X_n$ is $l$-extendible if and only if there exists a $\{0, 1, \ldots, l\}$-valued random variable $N$ such that

$$(4.1) \qquad P(X_1 = 1, \ldots, X_k = 1) = \frac{EN(N-1)\cdots(N-k+1)}{l(l-1)\cdots(l-k+1)}$$

for all $1 \leq k \leq n$.

To continue this program, we would need an analogue of (2.2) for $\{0, 1, \ldots, l\}$-valued random variables $N$. We are not aware of such a result. However, one can modify the approach that led to (2.2) to solve this problem. We will now do so, following the development in [13].

Fix integers $1 \leq n \leq l$. After a linear change of variables, our problem reduces to finding necessary and sufficient conditions on numbers $v_1, \ldots, v_n$ so that there exists a $\{0, 1, \ldots, l\}$-valued random variable $N$ satisfying

$$(4.2) \qquad v_k = EN^k \qquad \text{for } k = 1, \ldots, n.$$

The first step involves some convex analysis. Let $\mathcal{M}$ be the set of all nonnegative multiples of vectors $(1, v_1, \ldots, v_n)$, where $v_1, \ldots, v_n$ satisfy (4.2) for some $\{0, 1, \ldots, l\}$-valued random variable $N$, and let $\mathcal{P}$ be the set of polynomials $P(x)$ of degree at most $n$ that satisfy $P(i) \geq 0$ for $i = 0, \ldots, l$. Both $\mathcal{M}$ and $\mathcal{P}$ are closed convex cones. Clearly $P \in \mathcal{P}$ if and only if it has degree at most $n$ and $EP(N) \geq 0$ for all $\{0, \ldots, l\}$-valued $N$. Writing $P(x) = \sum_{i=0}^{n} c_i x^i$, we see that $P \in \mathcal{P}$ if and only if $\sum_{i=0}^{n} c_i v_i \geq 0$ for all $(v_0, \ldots, v_n) \in \mathcal{M}$. In other words, $\mathcal{P} = \mathcal{M}^*$, where $\mathcal{M}^*$ denotes the dual of $\mathcal{M}$. A basic result in convex analysis (see, e.g., Theorem 4.1 in [13]) then implies that $\mathcal{M} = \mathcal{P}^*$. This means that $(v_1, \ldots, v_n)$ can arise as in (4.2) if and only if $c_0 + \sum_{i=1}^{n} c_i v_i \geq 0$ for all $P(x) = \sum_{i=0}^{n} c_i x^i \in \mathcal{P}$.

It is sufficient in the last statement to consider only $P \in \mathcal{P}_e$, the extreme points of $\mathcal{P}$. To help understand the structure of $\mathcal{P}_e$, we state the following result.

PROPOSITION 4.1. *Suppose $P \in \mathcal{P}$. Then $P \in \mathcal{P}_e$ if and only if:*

(a) *$P$ has degree exactly $n$,*

*and*

(b) *all $n$ roots of $P$ are simple and are contained in $\{0, \ldots, l\}$.*

PROOF. Suppose $P \in \mathcal{P}_e$. Write $P(x) = P_1(x)P_2(x)$, where $P_1$ has only real roots, and $P_2$ has no real roots. Then $P_2$ is never zero, and we may assume that $P_2(x) > 0$ on $[0, l]$. If $P_2$ is not constant, we may write

$$P_2(x) = (P_2(x) + \varepsilon x)/2 + (P_2(x) - \varepsilon x)/2,$$

where the two polynomials on the right have degree at most that of $P_2(x)$ and are strictly positive on $[0, l]$ if $\varepsilon$ is sufficiently small. Then

$$P(x) = P_1(x)(P_2(x) + \varepsilon x)/2 + P_1(x)(P_2(x) - \varepsilon x)/2,$$

which violates the extremality of $P(x)$. Thus all roots of $P(x)$ are real.

Next, write $P(x) = (x - x_0)Q(x)$, where $x_0$ is one of the roots of $P(x)$. If $x_0 \notin \{0, \ldots, l\}$, then $(x - x_0 \pm \varepsilon)Q(x)$ are both in $\mathcal{P}$ for small $\varepsilon$, so that extremality is violated again. Therefore, all roots are in $\{0, \ldots, l\}$.

If $P(x)$ has degree less than $n$, then the representation

$$lP(x) = xP(x) + (l - x)P(x)$$

shows that $P$ is not extremal. Therefore, any extremal $P$ has degree $n$. One can rule out multiple roots at $i$ for $i = 0, 1, \ldots, l$ in a similar way, by writing

$$(x-i)^2 = \begin{cases} x(x-1) + x, & \text{if } i = 0, \\ (x-i)(x-i-1)/2 + (x-i)(x-i+1)/2, & \text{if } 1 \leq i \leq l-1, \\ (x-l)(x-l+1) + (l-x), & \text{if } i = l. \end{cases}$$

For the converse, suppose $P \in \mathcal{P}$ satisfies properties (a) and (b) in the statement of the proposition. Write $P = P_1 + P_2$ for $P_1, P_2 \in \mathcal{P}$. Then $P_1$ and $P_2$ have the same roots as $P$, and therefore must be positive multiples of $P$. This completes the proof of the proposition. □

Note that we have now established Theorem 1.16.

*Remarks.* With the above proof of Theorem 1.16, we can remove some of the mystery surrounding conditions (2.2) for the solution to the finite moment problem for a random variable $W$ taking values in $[0, 1]$. The proof in the continuous case is identical to that in the discrete case up to the statement of Proposition 4.1. However, in that case, the polynomials $P(x)$ in $\mathcal{P}_e$ have $n$ roots in $[0, 1]$, but all interior roots must have even multiplicity since $P(x) \geq 0$ on $[0, 1]$. If $n$ is even, for example, it follows that the roots at 0 or 1 (if any) must have multiplicities that are either both even or both odd. If they are both even, then $P(x) = [Q(x)]^2$ for some polynomial $Q(x)$, while if they are both odd, $P(x) = x(1-x)[Q(x)]^2$. Writing $Q(x) = \sum_i c'_i x^i$, replacing $x$ by $W$, squaring out these expressions and taking expected values, one sees that certain quadratic forms in the $c'_i s$ with coefficients that are $v_k$'s must be nonnegative definite. This (more or less) translates into the positivity of certain determinants, which turn out to be the ones in (2.2).

If $n = 2$, the polynomials appearing in Theorem 1.16 are just positive multiples of $x(l - x)$ and $(x - i)(x - i - 1)$ for $i = 0, \ldots, l - 1$. Therefore, the possibilities for $(v_1, v_2)$ are

$$\{(v_1, v_2) : v_2 \leq l v_1 \text{ and } v_2 \geq \phi(v_1)\},$$

where

$$\phi(x) = \max_{0 \leq i < l}[(2i + 1)x - i(i + 1)].$$

To see the analogy with (2.2), note that if $n = 2$, it becomes

$$v_2 > v_1^2 \quad \text{and} \quad v_2 < v_1,$$

and that

$$\lim_{l \to \infty} \phi(yl)/l^2 = y^2.$$

By (4.1), we see that a necessary and sufficient for $l$-extendibility of the Curie–Weiss Ising model with $n = 2$ is that there be a $\{0, 1, \ldots, l\}$-valued random variables $N$ so that

$$\frac{EN}{l} = \frac{a^2 + ab}{a^2 + 2ab + 1} \quad \text{and} \quad \frac{EN(N-1)}{l(l-1)} = \frac{a^2}{a^2 + 2ab + 1}.$$

By the previous development, this is equivalent to

$$\frac{la(b + la)}{a^2 + 2ab + 1} \geq \phi\left(\frac{la(a + b)}{a^2 + 2ab + 1}\right),$$

which is in turn equivalent to

$$la(b + la) \geq \max_{1 \leq i < l - 1}[(2i + 1)la(a + b) - i(i + 1)(a^2 + 2ab + 1)],$$

and hence to

(4.3) $$b \leq \min_{1 \leq i < l - 1}\left(\frac{(l - i)}{2i}a + \frac{i + 1}{2(l - i - 1)}a^{-1}\right).$$

(The inequalities corresponding to $i = 0$ and $i = l - 1$ are satisfied automatically.) Minimizing over continuous $i$ rather than discrete $i$ gives the following sufficient (and asymptotically necessary) condition for $l$-extendibility of the Curie–Weiss Ising model when $n = 2$:

$$b \leq 1 + \frac{(1 + a)^2}{2a(l - 1)}.$$

For $n = 3$, the polynomials appearing in Theorem 1.16 are $x(x - i)(x - i - 1)$ for $i = 1, \ldots, l - 1$ and $(l - x)(x - i)(x - i - 1)$ for $i = 0, \ldots, l - 2$. Therefore, the set of possible values of $(v_1, v_2, v_3)$ is given by

$$\{(v_1, v_2, v_3) : v_3 \geq \phi_1(v_1, v_2) \text{ and } 0 \leq v_3 \leq \phi_2(v_1, v_2)\},$$

where

$$\phi_1(x,y) = \max_{1 \leq i < l}[(2i+1)y - i(i+1)x]$$

and

$$\phi_2(x,y) = \min_{0 \leq i < l-1}[i(i+1)l - (2il + i^2 + i + l)x + (l + 2i + 1)y].$$

A necessary and sufficient condition for the $l$-extendibility of the Curie–Weiss Ising model with $n = 3$ is that there exist a $\{0, \ldots, l\}$-valued random variable $N$ so that

$$\frac{EN}{l} = \frac{a^3 + 2a^2b^2 + ab^2}{a^3 + 3a^2b^2 + 3ab^2 + 1},$$

$$\frac{EN(N-1)}{l(l-1)} = \frac{a^3 + a^2b^2}{a^3 + 3a^2b^2 + 3ab^2 + 1},$$

$$\frac{EN(N-1)(N-2)}{l(l-1)(l-2)} = \frac{a^3}{a^3 + 3a^2b^2 + 3ab^2 + 1}.$$

It then follows that this is equivalent to

$$(a^2 + b^2 + 2ab^2)i^2 - (2a^2l + 2ab^2l + b^2 - a^2)i + a(l-1)(2b^2 + al) \geq 0$$

for $1 \leq i \leq l-1$ and

$$(1 + 2ab^2 + a^2b^2)i^2$$
$$- (2ab^2l + 2a^2b^2l - 1 - 4ab^2 - 3a^3b^2)i + a^2b^2(l-1)(l-2) \geq 0$$

for $0 \leq i \leq l-2$. The values of the first quadratic at $i = 1, l-1$ are $a^2(l-1)(l-2)$ and $b^2(l-1)(l-2)$, respectively, and the values of the second quadratic at $i = 0, l-2$ are $a^2b^2(l-1)(l-2)$ and $(l-1)(l-2)$, respectively. Minimizing the quadratics over continuous values of $i$ as before, we see that a sufficient condition for $l$-extendibility of the Curie–Weiss Ising model for $n = 3$ is

$$b \leq 1 + \frac{\min[(a+b)^2, (1+ab)^2]}{2ab(l-2)}$$

provided that

$$l - 2 \geq \frac{\max[a^2 + b^2, 1 + a^2b^2]}{2ab^2}.$$

To check the assumption of Theorem 1.16 for larger values of $n$, we will need to appeal to part of the theory of quadratic forms. Recall from the remark following the statement of Theorem 1.16, that in the treatment of the continuous moment problem, one writes $P(x)$ as a perfect square, and uses the equivalence between the positive definiteness of a symmetric matrix and positivity of the principal minors of that matrix. In the discrete moment problem, the roots of $P$ are simple, so $P$ cannot be written as a perfect square. Nevertheless, the "interior" roots must appear as nearest neighbor pairs, so that they almost have even multiplicity. To quantify the difference that this makes, we need a quantitative version of the equivalence of positive definiteness and positivity of the principal minors. We turn to that next.

If $C = (c_{i,j})_{i,j \geq 0}$ is a matrix, we will use the following notation: $c_{-1,-1}^{-1} = 1$, and for $0 \leq k \leq i, j$,

$$c_{i,j}^k = \begin{vmatrix} c_{0,0} & \cdots & c_{0,k-1} & c_{0,j} \\ \cdots & \cdots & \cdots & \cdots \\ c_{k-1,0} & \cdots & c_{k-1,k-1} & c_{k-1,j} \\ c_{i,0} & \cdots & c_{i,k-1} & c_{i,j} \end{vmatrix}.$$

Note that this is the determinant of a $(k+1) \times (k+1)$ matrix because the indexes of the matrix $C$ begin at 0.

PROPOSITION 4.2. *If all of the principal minors $c_{k,k}^k$ of the symmetric matrix $C$ are nonzero, then*

$$(4.4) \qquad \sum_{i,j=0}^n c_{i,j} z_i z_j = \sum_{k=0}^n \frac{[\sum_{i=k}^n c_{i,k}^k z_i]^2}{c_{k,k}^k c_{k-1,k-1}^{k-1}}.$$

*Remarks.* (i) Note that this identity provides a simple proof of the standard fact referred to earlier that a quadratic form is positive definite if and only if the principal minors of the matrix of coefficients are all positive. The only if direction can be seen by perturbing the quadratic form a small amount such that all the principal minors are nonzero.

(ii) Equation (4.4) above is known as Jacobi's formula. Jacobi's original approach to it can be found in Chapter X, Section 3 of [9]—see equation (28) there. Nonetheless, we decided to include the proof here which seems to be different than the one in [9].

PROOF OF PROPOSITION 4.2. We begin with a special case of Sylvester's identity:

$$(4.5) \qquad c_{k,k}^k c_{i,j}^k - c_{i,k}^k c_{k,j}^k = c_{k-1,k-1}^{k-1} c_{i,j}^{k+1}, \qquad 0 \leq k \leq i, j,$$

where we have set $c_{i,j}^{k+1} = 0$ if $k = \min\{i,j\}$. This special case can be found at the bottom of page 586 of [1]. Proofs of the general identity can be found in that paper, as well as in [12].

Using (4.5), and the observation that the first sum below telescopes, we can write

$$c_{i,j} = \sum_{k=0}^{\min\{i,j\}} \left[ \frac{c_{i,j}^k}{c_{k-1,k-1}^{k-1}} - \frac{c_{i,j}^{k+1}}{c_{k,k}^k} \right] = \sum_{k=0}^{\min\{i,j\}} \frac{c_{i,k}^k c_{k,j}^k}{c_{k,k}^k c_{k-1,k-1}^{k-1}}.$$

Multiplying this identity by $z_i z_j$, summing, and then changing the order of summation gives

$$\sum_{i,j=0}^{n} c_{i,j} z_i z_j = \sum_{k=0}^{n} \sum_{i,j=k}^{n} \frac{c_{i,k}^{k} c_{k,j}^{k}}{c_{k,k}^{k} c_{k-1,k-1}^{k-1}} z_i z_j.$$

Finally, use the symmetry of $C$ to get (4.4). □

Now suppose $P$ is a polynomial of degree $n$ with $n$ simple roots in $\{0, \ldots, l\}$, and $P \geq 0$ on $\{0, \ldots, l\}$. Then the set of roots must be of the form $\{x_1, x_1 + 1\} \cup \cdots \cup \{x_p, x_p + 1\} \cup A$, where $A = \varnothing$ or $\{0, l\}$ if $n$ is even and $A = \{0\}$ or $\{l\}$ if $n$ is odd. To see this, suppose $P$ has consecutive simple roots at $\{x, x+1, \ldots, x+k-1\}$ and $P(x-1) > 0, P(x+k) > 0$. The sign of $P$ must change at $x, \ldots, x+k-1$, and therefore $k$ must be even. Now group these roots in pairs. For simplicity, for the moment we will take $n$ to be even. Similar results hold for odd $n$, but many of the formulas are a bit different. Unlike the "continuous" case described in the remark following the statement of Theorem 1.16, $(x - x_i)(x - x_i - 1)$ is not a perfect square, so in what follows, we will initially replace this product by $(x - x_i - \frac{1}{2})^2$.

Given $\{x_1, \ldots, x_p\}$, define $y_0 = 1$ and

$$y_q = (-1)^q \sum_{1 \leq i_1 < \cdots < i_q \leq p} (x_{i_1} + \tfrac{1}{2}) \cdots (x_{i_q} + \tfrac{1}{2})$$

for $1 \leq q \leq p$. Then we may write

(4.6)
$$\prod_{i=1}^{p} (x - x_i - \tfrac{1}{2})^2 = \left[\prod_{i=1}^{p} (x - x_i - \tfrac{1}{2})\right]^2$$
$$= \left(\sum_{i=0}^{p} x^{p-i} y_i\right)^2 = \sum_{i,j=0}^{p} x^{2p-i-j} y_i y_j.$$

Therefore, if we write

(4.7)
$$\prod_{i=1}^{p} (x - x_i - \tfrac{1}{2})^2 = \sum_{i=0}^{2p} c_i x^i,$$

it follows that for all $(v_0, \ldots, v_n) \in \mathcal{M}$

$$\sum_{i=0}^{2p} c_i v_i = \sum_{i,j=0}^{p} v_{2p-i-j} y_i y_j = \sum_{i,j=0}^{p} v_{i+j} y_{p-i} y_{p-j},$$

which is a quadratic form in the $y_{p-i}$'s. To apply Proposition 4.2, let $V$ be the matrix whose $(i,j)$ entry is $v_{i+j}$ for $i,j \geq 0$. If the necessary principal minors are nonzero, we then have

$$\text{(4.8)} \qquad \sum_{i=0}^{2p} c_i v_i = \sum_{k=0}^{p} \frac{[\sum_{i=k}^{p} v_{i,k}^k y_{p-i}]^2}{v_{k,k}^k v_{k-1,k-1}^{k-1}}.$$

This expression is used when $A = \varnothing$, in which case $n = 2p$. Similarly, if $A = \{0, l\}$, we have $p = (n-2)/2$ and we consider polynomials of the form

$$\text{(4.9)} \qquad x(l-x) \prod_{i=1}^{p} (x - x_i - \tfrac{1}{2})^2 = \sum_{i=0}^{2p+2} d_i x^i,$$

which yield, provided the necessary principal minors are nonzero,

$$\text{(4.10)} \qquad \sum_{i=0}^{2p+2} d_i v_i = \sum_{k=0}^{p} \frac{[\sum_{i=k}^{p} w_{i,k}^k y_{p-i}]^2}{w_{k,k}^k w_{k-1,k-1}^{k-1}},$$

where $W$ is the matrix whose $(i,j)$ entry is $w_{i+j+1}$ for $i, j \geq 0$, and $w_i = lv_i - v_{i+1}$ for $i \geq 1$.

But, of course, (4.7) and (4.9) are not the polynomials we must consider. So, we will use the identity

$$(x-u)(x-u-1) = (x - u - \tfrac{1}{2})^2 - \tfrac{1}{4}$$

to write the correct analogue of (4.7) as

$$\text{(4.11)} \qquad \begin{aligned} &\prod_{i=1}^{p} (x - x_i)(x - x_i - 1) \\ &= \sum_{q=0}^{p} (-\tfrac{1}{4})^q \sum_{1 \leq i_1 < \cdots < i_q \leq p} \prod_{j \neq i_1, \ldots, i_q} (x - x_j - \tfrac{1}{2})^2. \end{aligned}$$

For the correct analogue of (4.9), multiply both sides of (4.11) by $x(l-x)$.

Next we will illustrate the use of these expressions to check the assumptions of Theorem 1.16 for large $l$. It is this result that we will use in analyzing the Curie–Weiss Ising model. We will state it for even $n$; the case of odd $n$ is similar.

THEOREM 4.3. *Suppose $n = 2m$ is even and fixed, and $v_0(l) = 1, v_1(l), \ldots, v_n(l) > 0$ for $l \geq n$ are such that the corresponding quantities $v_{i,k}^k(l)$ and $w_{i,k}^k(l)$ are all positive and satisfy*

$$\frac{v_{k,k}^k(l)}{v_{k-1,k-1}^{k-1}(l)} = \gamma_k l^k + o(l^k), \qquad k = 0, \ldots, m,$$

$$\frac{w_{k,k}^k(l)}{w_{k-1,k-1}^{k-1}(l)} = \gamma_k' l^{k+2} + o(l^{k+2}), \qquad k = 0, \ldots, m-1,$$

$$\frac{v_{i,k}^k(l)}{v_{k,k}^k(l)} = \binom{i}{k}(\rho l)^{i-k} + o(l^{i-k}), \qquad 0 \leq k \leq i \leq n-k,$$

$$\frac{w_{i,k}^k(l)}{w_{k,k}^k(l)} = \binom{i}{k}(\rho l)^{i-k} + o(l^{i-k}), \qquad 0 \leq k \leq i \leq n-k-2,$$

as $l \to \infty$, where $0 < \rho < 1$ and $\gamma_k, \gamma_k' > 0$. Then for sufficiently large $l$, $v_1(l), \ldots, v_n(l)$ are the first $n$ moments of some $\{0, 1, \ldots, l\}$-valued random variable.

*Remarks.* (i) In view of the positivity of $v_{k,k}^k(l)$ and $w_{k,k}^k(l)$, it follows from criteria (2.2) that for such $l$ there is a random variable $N$ with values in $[0, l]$ with moments $v_0, \ldots, v_n$. The point of the above assumptions is to be guaranteed that $N$ can be taken to have values in $\{0, \ldots, l\}$.

(ii) If a family of finite moment sequences is $l$-extendible (in the sense of the conclusion of this theorem) for all large $l$, then (by passing to subsequences) there is at least one limiting distribution of an infinite exchangeable sequence. Theorem 4.3 is formulated for the situation in which this limiting distribution is the product measure of density $\rho$, because that is the case that arises in our analysis of the Curie–Weiss Ising model. One could presumably also use our technique to formulate analogous results for situations in which the limiting distribution is more general.

(iii) By the third display above with $k = 0$, we see that $v_i(l) = (\rho l)^i + o(l^i)$. Therefore the typical summand in the determinant $v_{i,k}^k$ is of order $(\rho l)^{k^2+i}$, and hence the expression in the first display above is potentially of order $(\rho l)^{2k}$. For it to be of order $l^k$ as is assumed here, there must be a lot of cancellation in the determinant. This is analogous to the fact that the variance of the sum of $m$ i.i.d. random variables is of order $m$ even though, without the cancellation that occurs, it would be of order $m^2$.

(iv) One can see from the proof that one actually only needs the third display to hold for $i \leq m$ and the fourth display to hold for $i \leq m - 1$.

PROOF OF THEOREM 4.3. We will apply Theorem 1.16. Consider a sequence $P_l(x)$ of polynomials of the form (4.7) in which $x_i = x_i(l) = \alpha_i l + o(l)$ for $i = 1, \ldots, p$, where $0 \leq \alpha_i \leq 1$ for each $i$. Then the corresponding $y_q$'s satisfy

$$y_q = (-l)^q \sum_{1 \leq i_1 < \cdots < i_q \leq p} \alpha_{i_1} \cdots \alpha_{i_q} + o(l^q).$$

Using the hypotheses of the theorem, it follows that the right-hand side of (4.8) equals

$$\sum_{k=0}^{p} \gamma_k (l^{2p-k} + o(l^{2p-k})) \left[ \sum_{1 \leq i_1 < \cdots < i_{p-k} \leq p} (\alpha_{i_1} - \rho) \cdots (\alpha_{i_{p-k}} - \rho) \right]^2.$$

Next we must account for the fact that the polynomial that arises in Theorem 1.16 is of form (4.11) rather than (4.7). If $c_i(l)$ is defined by

$$\prod_{I=1}^{m}(x - x_i(l))(x - x_i(l) - 1) = \sum_{i=0}^{n} c_i(l)x^i,$$

then

$$\sum_{i=0}^{n} c_i(l)v_i(l)$$
$$= (1 + o(1))$$
$$\times \sum_{B \subset \{1,\ldots,m\}} (-\tfrac{1}{4})^{m-|B|} \sum_{k=0}^{|B|} \gamma_k l^{2|B|-k} \left[ \sum_{D \subset B, |D|=|B|-k} \prod_{i \in D} (\alpha_i - \rho) \right]^2.$$

We need to check that this quantity is strictly positive for large $l$ for any choice of $\alpha_1, \ldots, \alpha_m \in [0, 1]$. The argument depends on how many of the $\alpha_i$'s are equal to $\rho$. For example, if $\alpha_i \neq \rho$ for all $1 \leq i \leq m$, then the term above corresponding to $B = \{1, \ldots, m\}$ and $k = 0$ is a positive multiple of $l^{2m}$, and all other terms are of smaller order. Suppose now that $\alpha_1 = \rho$ and $\alpha_i \neq \rho$ for all $2 \leq i \leq m$. Then the dominant term is of order $l^{2m-1}$, and corresponds to $B = \{1, \ldots, m\}$ and $k = 1$. More generally, suppose $\alpha_1 = \cdots = \alpha_j = \rho$ and $\alpha_i \neq \rho$ for all $i = j+1, \ldots, m$. Then the only $D$'s that can contribute to the expression satisfy $|D| \leq m - j$, since if $|D| > m - j$, one of the factors $\alpha_i - \rho$ must be zero. Therefore, for all nonzero summands above, $|B| - k \leq m - j$, and hence

(4.12) $\quad 2|B| - k = (|B| - k) + |B| \leq (m - j) + m \leq 2m - j.$

So, the largest power of $l$ that occurs in the above expression is $2m - j$. It can only occur if equality occurs in (4.12), that is, if $|B| = m$ and $k = j$. But in this case, the coefficient of $l^{2m-j}$ is $\gamma_j \prod_{i=j+1}^{m} (\alpha_i - \rho)^2$, which is strictly positive.

To complete the consideration of polynomials $P$ of the form (4.11) without assuming that $x_i(l)/l$ have limits, one passes to subsequences using compactness. The argument for polynomials of the form

$$x(l - x) \prod_{i=1}^{m-1} (x - x_i)(x - x_i - 1),$$

is similar, using the assumptions on the $w$'s rather than the $v$'s. We have now verified that if $P_l(x) = \sum_{i=0}^{n} c_i(l)x^i$ is a polynomial for each $l \geq n$ that has $n$ simple roots in $\{0, \ldots, l\}$ and is nonnegative on $\{0, \ldots, l\}$, then for sufficiently large $l$, $\sum_{i=0}^{n} c_i(l)v_i(l) \geq 0$. It follows from Theorem 1.16 that for such $l$, $v_1(l), \ldots, v_n(l)$ are the first $n$ moments of some $\{0, \ldots, l\}$-valued random variable. $\square$

We are now ready for the

PROOF OF THEOREM 1.15. For the Curie–Weiss Ising model, we have

$$P(X_1 = 1, \ldots, X_k = 1) = \sum_{j=k}^{n} \binom{n-k}{j-k} u_j, \qquad 1 \leq k \leq n,$$

where $u_j = a^j b^{j(n-j)}/s_n$ and

$$s_n = \sum_{j=0}^{n} \binom{n}{j} a^j b^{j(n-j)}.$$

With this notation, (4.1) becomes

(4.13)
$$\begin{aligned} &EN(N-1)\cdots(N-k+1) \\ &= l(l-1)\cdots(l-k+1) \sum_{j=k}^{n} \binom{n-k}{j-k} u_j, \qquad 1 \leq k \leq n. \end{aligned}$$

Expanding the product on the left-hand side of (4.13) and writing $v_k(l) = EN^k$, we can solve the resulting equations for these quantities. The issue is whether $v_1(l), \ldots, v_n(l)$ are the first $n$ moments of a $\{0, \ldots, l\}$-valued random variable.

Defining $w_k(l) = lv_k(l) - v_{k+1}(l)$, we will apply Theorem 4.3 to prove part (a), and therefore assume at this point that $n$ is even. To verify the assumptions of Theorem 4.3, and also for the easy proof of part (b), we will need the following asymptotic statements: As $l \to \infty$,

$$(4.14) \quad v_{k,k}^k(l) \sim \left(\prod_{j=0}^{k} j!\right) [\delta\rho(1-\rho)l]^{k(k+1)/2},$$

$$(4.15) \quad w_{k,k}^k(l) \sim \left(\prod_{j=0}^{k} j!\right) \delta^{k(k+1)/2} [\rho(1-\rho)]^{(k+1)(k+2)/2} l^{(k+1)(k+4)/2},$$

$$(4.16) \quad v_{i,k}^k(l) \sim \binom{i}{k} (\rho l)^{i-k} v_{k,k}^k(l) \quad \text{and} \quad w_{i,k}^k(l) \sim \binom{i}{k} (\rho l)^{i-k} w_{k,k}^k(l),$$

where $\delta = 1 - 2c\rho(1-\rho)$.

The hypothesis of part (a) of the theorem gives $\delta > 0$, which is what we need to apply Theorem 4.3. Part (b) of the theorem follows immediately from (4.14) with $k = 1$, since

$$v_{1,1}^1 = \begin{vmatrix} v_0 & v_1 \\ v_1 & v_2 \end{vmatrix} = v_2 - v_1^2,$$

which is nonnegative if $v_1, v_2$ are the first two moments of any random variable.

To check (4.14), (4.15) and (4.16), we need to solve (4.13) explicitly. In order to do so, let $G = (g_{i,j})_{0 \leq i,j \leq n}$ and its inverse $H = (h_{i,j})_{0 \leq i,j \leq n}$ be defined by

$$G \begin{pmatrix} 1 \\ l \\ l^2 \\ \vdots \\ l^n \end{pmatrix} = \begin{pmatrix} 1 \\ l \\ l(l-1) \\ \vdots \\ l(l-1)\cdots(l-n+1) \end{pmatrix}$$

and

$$H \begin{pmatrix} 1 \\ l \\ l(l-1) \\ \vdots \\ l(l-1)\cdots(l-n+1) \end{pmatrix} = \begin{pmatrix} 1 \\ l \\ l^2 \\ \vdots \\ l^n \end{pmatrix}$$

for every $l$. These are the lower triangular matrices

$$G = \begin{pmatrix} 1 & 0 & \cdots & 0 \\ 0 & 1 & \cdots & 0 \\ 0 & -1 & \cdots & 0 \\ 0 & 2 & \cdots & 0 \\ \vdots & \vdots & \ddots & \vdots \\ 0 & (-1)^{n-1}(n-1)! & \cdots & 1 \end{pmatrix}$$

and

$$H = \begin{pmatrix} 1 & 0 & 0 & \cdots & 0 \\ 0 & 1 & 0 & \cdots & 0 \\ 0 & 1 & 1 & \cdots & 0 \\ 0 & 1 & 3 & \cdots & 0 \\ \vdots & \vdots & \vdots & \ddots & \vdots \\ 0 & 1 & 2^{n-1}-1 & \cdots & 1 \end{pmatrix}.$$

Entries in these matrices other than those given above are rather complicated. However, they are completely determined by equating coefficients of powers of $l$ in the defining relations above; we will later give recursive expressions for them. To simplify the following expressions, we will often suppress the dependence on $l$, and write $v_k = v_k(l)$. We will also suppress the limits of the following sums, relying on the usual convention that $\binom{m}{k} = 0$ except when $0 \leq k \leq m$. Then the solution of (4.13) is given by

$$(4.17) \qquad v_m = \sum_{i,j,k} h_{m,k} g_{k,j} l^j \binom{n-k}{i-k} u_i.$$

To motivate the next step, recall that we are trying to prove that for large $l$, the $v_m$'s are the moments of a random variable $N$ that has a distribution close to $B(l, \rho)$. If this were the case, then $E(N - l\rho)^p$ should be of order $l^{p/2}$ rather than $l^p$. We will now check this without assuming that the $v_m$'s are moments at all. Recalling that $\rho = a/(1+a)$, it is natural to consider the following, which we rewrite using (4.17):

$$\sum_m \binom{p}{m}(1+a)^m(-la)^{p-m}v_m$$

$$= \sum_{i,j,k,m} \binom{p}{m}(1+a)^m(-la)^{p-m}h_{m,k}g_{k,j}l^j\binom{n-k}{i-k}u_i.$$

Next write $b = 1 + (c/l)$ and use the binomial expansion for powers of $1+(c/l)$ to write

(4.18)
$$s_n \sum_m \binom{p}{m}(1+a)^m(-la)^{p-m}v_m$$
$$= \sum_{i,j,k,m,q} \binom{p}{m}(1+a)^m(-la)^{p-m}h_{m,k}g_{k,j}l^j$$
$$\times \binom{n-k}{i-k}a^i\binom{i(n-i)}{q}(c/l)^q$$
$$= \sum_{r,s,t:r,s\geq 0, r+s+t\leq p} l^t a^{p-s}(1+a)^s c^{p-r-s-t}(-1)^{p-s}C_{r,s,t}(a),$$

where

$$C_{r,s,t}(a) = \sum_{k=0}^p \sum_{i=k}^n a^{i-k}(1+a)^k \binom{p}{k+s}(-1)^k$$
$$\times \binom{i(n-i)}{p-r-s-t}\binom{n-k}{i-k}h_{s+k,k}g_{k,k-r}.$$

In the final step, we have let $r = k - j, s = m - k$ and $t = j - q + p - m$. Shortly, we will show that

(4.19) $C_{r,s,t}(a) = \begin{cases} 0, & \text{if } t > p/2, \\ (1\cdot 3\cdots(p-1))(-1)^t\binom{t}{r+s}\binom{r+s}{r} \\ \quad \times 2^{t-r-s}(1+a)^{n-p+2r+2s}, & \text{if } t = p/2. \end{cases}$

But first, we will use this to complete the verification of (4.14).

Since $s_n \to (1+a)^n$ as $l \to \infty$, it follows from (4.19) that

(4.20) $\sum_m \binom{p}{m}(-\rho l)^{p-m}v_m \begin{cases} = O(l^{(p-1)/2}), \\ \quad \text{if } p \text{ is odd}, \\ \sim [1\cdot 3\cdots(p-1)][\delta\rho(1-\rho)l]^{p/2}, \\ \quad \text{if } p \text{ is even}. \end{cases}$

Applying row and column operations similar to what was done in the alternative proof of Theorem 1.1, we see that $v_{k,k}^k(l)$ can be written as

$$\begin{vmatrix} v_0 & v_1 - \rho l v_0 & \cdots & \sum_m \binom{k}{m}(-\rho l)^{k-m} v_m \\ v_1 - \rho l v_0 & v_2 - 2\rho l v_1 + \rho^2 l^2 v_0 & \cdots & \sum_m \binom{k+1}{m}(-\rho l)^{k+1-m} v_m \\ \vdots & \vdots & \ddots & \vdots \\ \sum_m \binom{k}{m}(-\rho l)^{k-m} v_m & \sum_m \binom{k+1}{m}(-\rho l)^{k+1-m} v_m & \cdots & \sum_m \binom{2k}{m}(-\rho l)^{2k-m} v_m \end{vmatrix}.$$

Using (4.20), we see that

$$v_{k,k}^k(l) \sim [\delta \rho (1-\rho) l]^{1+2+\cdots+k} \begin{vmatrix} 1 & EZ & EZ^2 & \cdots & EZ^k \\ EZ & EZ^2 & EZ^3 & \cdots & EZ^{k+1} \\ \vdots & \vdots & \vdots & \ddots & \vdots \\ EZ^k & EZ^{k+1} & EZ^{k+2} & \cdots & EZ^{2k} \end{vmatrix},$$

where $Z$ is a standard normal random variable. By Corollary 4C in [17], this last determinant is $1!2!\cdots k!$. It can also be deduced from results on the $\beta = 2$ Gaussian ensemble—see Chapters 4 and 17 of [18].

It remains to prove (4.19). First, we need some information about the matrices $G$ and $H$. By the definition of $G$, for example,

$$\sum_{j=0}^{i} g_{i,j} l^j = l(l-1)\cdots(l-i+1).$$

Therefore,

$$\sum_{j=0}^{i+1} g_{i+1,j} l^j = (l-i) \sum_{j=0}^{i} g_{i,j} l^j = \sum_{j=1}^{i+1} g_{i,j-1} l^j - i \sum_{j=0}^{i} g_{i,j} l^j.$$

Equating coefficients of powers of $l$ gives $g_{i+1,j} = g_{i,j-1} - i g_{i,j}$ and then solving the resulting recursion leads to

$$g_{i+1,j} = -\sum_{k=1}^{j} g_{i-j+k,k}(i-j+k).$$

Similarly, $h_{i+1,j} = h_{i,j-1} + j h_{i,j}$ and

$$h_{i+1,j} = \sum_{k=1}^{j} k h_{i-j+k,k}.$$

It follows by induction that

(4.21)
$$g_{k,k-r} = \frac{(-1)^r}{2^r r!} \times \text{a monic polynomial in } k \text{ of degree } 2r,$$

$$h_{s+k,k} = \frac{1}{2^s s!} \times \text{a monic polynomial in } k \text{ of degree } 2s.$$

The next step is to prove by induction on $j$ that if $P$ is any monic polynomial in $i$ of degree $j$, then

$$\sum_i P(i) a^{i-k} \binom{n-k}{i-k} = Q(k)(1+a)^{n-k-j} \qquad (4.22)$$

for some monic polynomial $Q$ in $k$ of degree $j$. This is clearly true for $j = 0$. For the induction step, write

$$\sum_i i^j a^{i-k} \binom{n-k}{i-k}$$

$$= \sum_i i^{j-1}[(i-k) + k] a^{i-k} \binom{n-k}{i-k}$$

$$= (n-k)a \sum_i i^{j-1} a^{i-k-1} \binom{n-k-1}{i-k-1} + k \sum_i i^{j-1} a^{i-k} \binom{n-k}{i-k}$$

$$= (n-k)a Q_1(k+1)(1+a)^{n-(k+1)-(j-1)} + k Q_2(k)(1+a)^{n-k-(j-1)}$$

$$= Q(k)(1+a)^{n-k-j},$$

where $Q_1(x)$ and $Q_2(x)$ are monic polynomials of degree $j - 1$ (by the inductive hypothesis) and

$$Q(k) = (n-k)a Q_1(k+1) + k(1+a) Q_2(k)$$

is a monic polynomial of degree $j$.

Now start with

$$\sum_k \binom{p}{k}(-1)^k x^k = (1-x)^p,$$

differentiate $j$ times with respect to $x$ and set $x = 1$. The result is that

$$\sum_k \binom{p}{k}(-1)^k k(k-1)\cdots(k-j+1) = \begin{cases} 0, & \text{if } j < p, \\ (-1)^j p!, & \text{if } j = p. \end{cases}$$

It follows that

$$\sum_k \binom{p}{k}(-1)^k k^j = \begin{cases} 0, & \text{if } j < p, \\ (-1)^j p!, & \text{if } j = p. \end{cases} \qquad (4.23)$$

Since $\binom{i(n-i)}{p-r-s-t}$ is a polynomial in $i$ of degree $2(p - r - s - t)$ and leading coefficient $(-1)^{p-r-s-t}/(p-r-s-t)!$, (4.22) implies that

$$\sum_i a^{i-k} \binom{i(n-i)}{p-r-s-t} \binom{n-k}{i-k}$$

$$= \frac{(-1)^{p-r-s-t}(1+a)^{n-k-2(p-r-s-t)}}{(p-r-s-t)!} Q(k),$$

where $Q(k)$ is a monic polynomial of degree $2(p - r - s - t)$. Therefore,

$$C_{r,s,t}(a) = \frac{(-1)^{p-r-s-t}(1+a)^{n-2(p-r-s-t)}}{(p-r-s-t)!}$$
$$\times \sum_k \binom{p}{k+s}(-1)^k Q(k) h_{s+k,k} g_{k,k-r}$$
$$= \frac{(-1)^{p-s-t}(1+a)^{n-2(p-r-s-t)}}{2^{r+s} r! s! (p-r-s-t)!} \sum_k \binom{p}{k+s}(-1)^k \times Q^*(k)$$
$$= \frac{(-1)^{p-t}(1+a)^{n-2(p-r-s-t)}}{2^{r+s} r! s! (p-r-s-t)!} \begin{cases} 0, & \text{if } 2(p-t) < p, \\ (2(p-t))!, & \text{if } 2(p-t) = p. \end{cases}$$

Here $Q^*(k)$ is a monic polynomial of degree $2(p-t)$. We have used (4.21) in the middle equality and (4.23) in the final one.

This completes the proof of (4.14). The proof of (4.15) is similar. In particular, we now know that for large $l$, $v_{k,k}^k(l)$ and $w_{k,k}^k(l)$ are strictly positive. It follows from criteria (2.2) that for such $l$ there is a random variable $N$ with values in $[0,l]$ with moments $v_0, \ldots, v_n$. We still do not know it can be taken to have values in $\{0, \ldots, l\}$, and for that we still need to check (4.16). Note that the expression for $v_{k,k}^k(l)$ that follows (4.20) can then be written as

$$\begin{vmatrix} 1 & E(N-\rho l) & \cdots & E(N-\rho l)^k \\ E(N-\rho l) & E(N-\rho l)^2 & \cdots & E(N-\rho l)^{k+1} \\ \vdots & \vdots & \ddots & \vdots \\ E(N-\rho l)^k & E(N-\rho l)^{k+1} & \cdots & E(N-\rho l)^{2k} \end{vmatrix}.$$

Similarly, row and column operations can be used to write $v_{i,k}^k(l)$ as

$$(4.24) \quad \begin{vmatrix} 1 & E(N-\rho l) & \cdots & E(N-\rho l)^k \\ E(N-\rho l) & E(N-\rho l)^2 & \cdots & E(N-\rho l)^{k+1} \\ \vdots & \vdots & \ddots & \vdots \\ E(N-\rho l)^{k-1} & E(N-\rho l)^k & \cdots & E(N-\rho l)^{2k-1} \\ E(N-\rho l)^k f_k(\rho l, N) & E(N-\rho l)^{k+1} f_k(\rho l, N) & \cdots & E(N-\rho l)^{2k} f_k(\rho l, N) \end{vmatrix},$$

where

$$f_k(x,y) = \sum_{j=0}^{i-1} j(j-1)\cdots(j-k+2) x^{j-k+1} y^{i-j-1}.$$

Note that while it might appear that $E(N-\rho l)^k f_k(\rho l, N)$, for example, is a linear combination of all the moments of $N$ up to order $i$, and therefore would not be obtainable as a linear combination of the moments of order $0, 1, \ldots, k-1$ and $i$ that originally appeared in the first column of the matrix, it really is a linear combination of these latter moments. This can be seen by writing

$$f_k(x,y) = \frac{d^{k-1}}{dx^{k-1}} \frac{y^i - x^i}{y - x}.$$

To complete the proof of the first part of (4.16), note that by (4.20), $N/l \to \rho$ in probability, and therefore,

$$\frac{f_k(\rho l, N)}{(\rho l)^{i-k}} \to \sum_{j=0}^{i-1} j(j-1)\cdots(j-k+2) = \binom{i}{k}$$

in probability. It follows that asymptotically, one can factor out a term

$$\binom{i}{k}(\rho l)^{i-k}$$

from the last row of the above matrix. This gives the first part of (4.16). The proof of the last part is similar. $\square$

*Remarks.* (i) In case (a) of Theorem 1.15, we can consider for large $l$ a random variable $N_l$ that corresponds to the extension of length $l$ of the Curie–Weiss Ising model. (This is presumably not unique.) By (4.20), any weak limit of

$$\frac{N_l - \rho l}{\sqrt{\delta \rho (1 - \rho) l}}$$

as $l \to \infty$ has the same first $n$ moments as a standard normal random variable.

(ii) Diaconis and Freedman [5] proved that if the exchangeable measure $\mu$ on $\{0,1\}^n$ is $l$-extendible, then the total variation distance between $\mu$ and the closest mixture of homogeneous product measures is at most $4n/l$. Combining this statement with Theorem 1.15 yields the following conclusion: If $\mu_l$ is the distribution of the Curie–Weiss Ising model on $\{0,1\}^n$ with $a = \rho/(1-\rho)$ and $b = 1 + (c/l)$, where $c < 1/(2\rho(1-\rho))$, then

$$\limsup_{l \to \infty} l \inf_{\gamma} \left\| \mu_l - \int_0^1 \nu_\tau \gamma(d\tau) \right\|_{\mathrm{TV}} \leq 4n,$$

where $\nu_\tau$ is the homogeneous product measure with density $\tau$, and the infimum is over probability measures $\gamma$ on $[0,1]$.

One could try to take our analysis further at the critical point $c = 1/(2\rho(1-\rho))$ (with $\rho \neq 0, 1$). To see what can happen, let

(4.25) $$b = 1 + \frac{1}{2\rho(1-\rho)l} + \frac{d}{l^2}.$$

Theorem 4.3 no longer applies, but for small $n$, one can apply Theorem 1.16 directly, as we did earlier in this section. Here are some results for $n = 2$, where (4.3) is being used, that illustrate the complexity of the answer: Suppose $\rho$ is rational, and write $\rho = \frac{j}{k}$. Take integers $m \in [0, k/2]$ and $q \in [0, k)$ so that $jq \equiv m \pmod{k}$ or $jq \equiv k - m \pmod{k}$. Then the critical value for $l$-extendibility for all large $l \equiv q \pmod{k}$ is

$$d_c = \frac{\rho(1-\rho) + (m/k)^2}{2\rho^2(1-\rho)^2}.$$

Perhaps surprisingly, if $n \geq 4$ and $\rho \neq \frac{1}{2}$, there is no value of $d$ for which $l$-extendibility holds for all large $l$. To see this, one can solve (4.13) for $k \leq 4$ and check that $N_l$, if it exists, must satisfy the following, as $l \to \infty$:

$$E(N_l - \rho l) \to \frac{n-1}{2}(1 - 2\rho),$$

$$E(N_l - \rho l)^2 \to \frac{(n^2 - 4n + 6) - 4(n^2 - 5n + 6)\rho(1-\rho)}{4} - 2d\rho^2(1-\rho)^2,$$

$$E(N_l - \rho l)^3 \sim \rho(1-\rho)(1-2\rho)l,$$

$$E(N_l - \rho l)^4 \sim 2\rho(1-\rho)[n(1-2\rho)^2 + 13\rho(1-\rho) - 3]l.$$

If $\rho \neq \frac{1}{2}$, this violates the Schwarz inequality: $(E(N - \rho l)^3)^2 \leq E(N - \rho l)^2 E(N - \rho l)^4$. If $\rho = \frac{1}{2}$, one must compute higher moments to draw the same conclusion (for $n \geq 6$):

$$E(N_l - \rho l)^2 \to \frac{2n-d}{8},$$

$$E(N_l - \rho l)^4 \sim \tfrac{1}{8}l,$$

$$E(N_l - \rho l)^6 \sim \frac{30n - 15d - 56}{64}l.$$

Presumably, this means that the power of $l$ that is used in the correction in (4.25) is not necessarily 2, and may depend on $\rho$ and/or $n$. We have not investigated this further.

**5. A formula for finite extensions.** In this section, we prove Propositions 1.17, 1.18, 1.19 and 1.20 as well as Theorem 1.21.

PROOF OF PROPOSITION 1.17. Fix $n$, $J > 0$, $h \geq 0$ and $l > n$. It is easy to check that for $k = 0, 1, \ldots, n$, the probability that there are $k$ 1's in the Curie–Weiss Ising model with parameters $n$, $-J$ and $h$ is given by

$$\frac{1}{Z_n}\binom{n}{k}\exp\{-J(2k-n)^2/2 + h(2k-n)\}$$

$$= \frac{1}{Z_n}\binom{n}{k}\int_{-\infty}^{\infty} e^{(ix+h)(2k-n)} f_J(x)\,dx.$$

Now, define for each $j \in \{0, \ldots, l\}$
$$\tilde{Q}(j) := \frac{1}{Z_n} \binom{l}{j} \int_{-\infty}^{\infty} e^{(ix+h)(2j-l)} (e^{(ix+h)} + e^{-(ix+h)})^{n-l} f_J(x) \, dx.$$

For $h \neq 0$, there is no singularity and so the integral is well defined. It is not hard to see that $\tilde{Q}(j)$ is real but our later computation will verify this. It is also straightforward to verify the claim concerning (1.12) for $h > 0$ with $Q$ replaced by $\tilde{Q}$. (It is also true that the $\tilde{Q}$ corresponding to $l' > l$ has the $\tilde{Q}$ corresponding to $l$ as hypergeometric projection in the obvious sense. Although not so interesting, the above makes sense and is also correct when $l \leq n$.) The proof will be complete if we show that $\tilde{Q} = Q$ when $h > 0$. The $h = 0$ case of (1.12) will then follow from continuity.

To see that $\tilde{Q} = Q$ when $h > 0$, we now let $M := l - k > 0$ and $u := 2j - l \in \{-l, -l+2, \ldots, l-2, l\}$. If $h > 0$, we have

$$\int_{-\infty}^{\infty} e^{(ix+h)u} (e^{(ix+h)} + e^{-(ix+h)})^{-M} f_J(x) \, dx$$

$$= \int_{-\infty}^{\infty} e^{(ix+h)(u-M)} (1 + e^{-2(ix+h)})^{-M} f_J(x) \, dx$$

$$= \int_{-\infty}^{\infty} e^{(ix+h)(u-M)} \sum_{m=0}^{\infty} \binom{-M}{m} e^{-2m(ix+h)} f_J(x) \, dx$$

$$= \sum_{m=0}^{\infty} \binom{-M}{m} e^{-h(2m+M-u)} \int_{-\infty}^{\infty} e^{-ix(2m+M-u)} f_J(x) \, dx$$

$$= \sum_{m=0}^{\infty} \binom{-M}{m} e^{-h(2m+M-u)} e^{-(J/2)(2m+M-u)^2}$$

$$= e^{h^2/(2J)} \sum_{m=0}^{\infty} \binom{-M}{m} e^{-J/2(2m+M-u+h/J)^2},$$

as desired. □

PROOF OF PROPOSITION 1.18. First, assume that $n$ is odd and $h = 0$. In this case, we will apply Proposition 1.17 by simply verifying that the expression given there is nonnegative. The series we need to consider is

$$\sum_{m=0}^{\infty} (-1)^m e^{-(J/2)(2m+1+u)^2}.$$

When $u \geq 0$, this is nonnegative since it is an alternating series with terms whose absolute values are decreasing. If $u = -2j$, where $j = 1, 2, \ldots$, write this series as

$$\sum_{m=0}^{2j-1} (-1)^m e^{-(J/2)(2m+1+u)^2} + \sum_{m=2j}^{\infty} (-1)^m e^{-(J/2)(2m+1+u)^2}.$$

The second sum is nonnegative again because it is an alternating series with terms whose absolute values are decreasing. The first sum is zero, since the summands for $m$ and $2j - m - 1$ cancel.

For the converse, we use our $a, b$ parameterization. Fix $n$, $a$, $b$ and $l$. If the Curie–Weiss Ising model with parameters $a$ and $b$ on $\{0, 1\}^n$ is $l$-extendible, then there must exist a random variable $N$ taking values in $\{0, \ldots, l\}$ satisfying (4.13). Using this equation, we can compute the variance of $N$ which then turns out to be

$$(\ell/s_n^2) \sum_{i,j} \left[ \ell \binom{n-2}{j-2} \binom{n}{i} + \binom{n-2}{j-1} \binom{n}{i} \right.$$
$$\left. - \ell \binom{n-1}{i-1} \binom{n-1}{j-1} \right] a^{i+j} b^{i(n-i)+j(n-j)}.$$

It is clear that all cases fall into one of the following three cases.

Case (1): $n$ is even, $\ell = n+1$ and $a$ is arbitrary. In this case, the dominant term (as $b$ gets large) is $i = j = n/2$ which can be seen to have a negative coefficient. Hence for large $b$ the above is negative and the extension does not exist.

Case (2): $n$ is odd, $\ell = n + 2$ and $a$ is arbitrary. In this case, there are four dominant terms (as $b$ gets large) corresponding to $i, j \in \{(n-1)/2, (n+1)/2\}$. An easy computation shows that in this case the sum of the coefficients of these four terms is negative and hence for large $b$ the above is negative and the extension does not exist.

Case (3): $n$ is odd, $\ell = n + 1$ and $a \neq 1$. In this case, there are again the same four dominant terms (as $b$ gets large) as in the previous case and an easy computation again shows that in this case the sum of the coefficients of these four terms is negative and so, as before, the extension does not exist. □

PROOF OF PROPOSITION 1.19. Again, the series we need to consider is

$$\sum_{m=0}^{\infty} (-1)^m e^{-(J/2)(2m+1+u)^2}.$$

If $u \geq -1$, then this is an alternating series with terms whose absolute values are decreasing and hence is nonnegative. Otherwise, write $u = -2j - 1$ with $j = 1, 2, 3, \ldots$ where there are only finitely many $j$'s here. Then the above sum becomes

$$\sum_{m=0}^{\infty} (-1)^m e^{-2J(m-j)^2}.$$

Break the sum into $m \leq 2j + 1$ and $m \geq 2j + 2$. The second sum is fine as before. For the first sum, expand the exponential in powers of $J$. The constant term is 0 because there are an even number of summands. The coefficient of $J$ in the expansion is

$$2 \sum_{m=0}^{2j+1} (-1)^{m+1}(m-j)^2 = 2(j+1) > 0.$$

This term dominates for small $J$ and hence the sum is positive. We now apply Proposition 1.17. □

PROOF OF PROPOSITION 1.20. Fix $n$, $l$ and $p \in (0,1)$. The measures in $E(l,l)$ are $l$-dimensional and correspond to a simplex $A$ in $R^l$. Similarly, the measures in $E(n,n)$ are $n$-dimensional and correspond to a simplex $B$ in $R^n$. The hypergeometric projection corresponds to a linear mapping $f$ from $A$ to $B$ whose image is a set $C$ in $R^n$ corresponding exactly to $E(n,l)$. Clearly as a map from $A$ to $C$, $f$ has full rank. Consider the point $a$ in $A$ corresponding to the process $Y = (Y_1, \ldots, Y_l)$ which is i.i.d. 0, 1 valued with $P(Y_1 = 1) = p$. Clearly $a$ is an interior point of $A$. Since $f$ has full rank, $f(a)$ is an interior point of $C$. However, $f(a)$ corresponds exactly to $X = (X_1, \ldots, X_n)$. This proves the claim.

For the last statement, fix $p = 1$. We take $n = 4$ and $l = 5$. Consider the measure on $\{0,1\}^4$ which is $(1-\varepsilon)m_1 + \varepsilon m_2$ where $m_1$ is product measure with $p = 1$ and $m_2$ is uniform distribution on configurations with exactly 2 1's. It is easy to see that for any $\varepsilon > 0$, this measure, while in $E(4,4)$, is not 5-extendible. □

REMARK. There is an alternative way to prove the above result. When one extends the product measure from $n$ to $l$ sites, the resulting random variable $N$ satisfying (4.1) is binomially distributed. By Proposition 4.1, we therefore have that for every polynomial $P \in \mathcal{P}_e$, $EP(N) > 0$. Since $|\mathcal{P}_e| < \infty$, it follows by Theorem 1.16, that if the finite sequence $v_0, \ldots, v_n$ is close to these binomial moments, then they are also the moments of some $N'$ of the desired form. This, together with (4.13), completes the alternate proof. In fact, this proof shows that whenever we have a process $\{X_1, \ldots, X_n\}$ in $E(n,l)$ which has some "representing" $N$ satisfying (4.1) having at least $n+1$ points in its support, then small perturbations of $\{X_1, \ldots, X_n\}$ which are in $E(n,n)$ are also in $E(n,l)$.

We finally now move to the proof of Theorem 1.21.

PROOF OF THEOREM 1.21. Fix $n, c, h$ and $l$. Letting $J = c/l$, we have, using (1.11), $Q(j)$ defined for each $j \in \{0, \ldots, l\}$. We want to show that for $h > h^*(c)$, we have that for large $l$, $Q(j)$ is nonnegative for all $j \in \{0, \ldots, l\}$.

Since $h^*(c) > 0$, in view of the proof of Proposition 1.17, we need to show that for $j \in \{0, \ldots, l\}$,

$$\tilde{Q}(j) := \int_{-\infty}^{\infty} e^{(ix+h)(2j-l)} (e^{(ix+h)} + e^{-(ix+h)})^{n-l} f_J(x)\, dx \geq 0,$$

where we recall that $f_J$ is the density function for a normal random variable with mean 0 and variance $J$.

Since we do not care about positive multiplicative factors, we will use the notation $A \cong B$ if the $A$ and $B$ only differ by a positive multiplicative factor. Letting $2j - l = vl$ with $v \in [-1, 1]$, a simple change of variables shows that

$$(5.1) \qquad \tilde{Q}(j) \cong \int_{-\infty}^{\infty} e^{ixv\sqrt{l}} (e^{h+ix/\sqrt{l}} + e^{-h-ix/\sqrt{l}})^{n-l} f_c(x)\, dx.$$

Thinking of $x$ as *complex*, the integrand on the right-hand side of (5.1) has isolated poles of order $l - n$ at the points $x = (ih + (2r+1)\pi/2)\sqrt{l}$ for $r \in \mathbb{Z}$ and is otherwise analytic. One can then readily deduce from Cauchy's theorem that if $h, \chi > 0$, then the integrand in (5.1) is unchanged if we integrate over $\mathbb{R} + i(h - \chi)\sqrt{l}$ instead. This leads to

$$(5.2) \quad \begin{aligned} \tilde{Q}(j) &\cong \int_{-\infty}^{\infty} e^{ix(v-(h-\chi)/c)\sqrt{l}} (e^{\chi+ix/\sqrt{l}} + e^{-\chi-ix/\sqrt{l}})^{n-l} \varphi_c(x)\, dx \\ &\cong \int_{-\infty}^{\infty} (e^{\chi+ix/\sqrt{l}} + e^{-\chi-ix/\sqrt{l}})^n \\ &\qquad \times (e^{\chi+(ix/\sqrt{l})(1-v+(h-\chi)/c)} \\ &\qquad\quad + e^{-\chi-(ix/\sqrt{l})(1+v-(h-\chi)/c)})^{-l} \varphi_c(x)\, dx. \end{aligned}$$

Now, assume $\chi = \chi(c, h, v)$ is a solution of the equation

$$(5.3) \qquad \xi - c\tanh\xi = h - cv, \qquad \xi > 0.$$

Denote

$$p = p(c, h, v) := \frac{e^{\chi}}{e^{\chi} + e^{-\chi}}, \qquad q = q(c, h, v) := \frac{e^{-\chi}}{e^{\chi} + e^{-\chi}},$$

and observe that

$$p + q = 1, \qquad p - q = \tanh\chi \quad \text{and} \quad 4pq = \frac{1}{(\cosh\chi)^2}.$$

We then readily obtain

$$(5.4) \quad \tilde{Q}(j) \cong \int_{-\infty}^{\infty} (pe^{ix/\sqrt{l}} + qe^{-ix/\sqrt{l}})^n (pe^{i2qx/\sqrt{l}} + qe^{-i2px/\sqrt{l}})^{-l} \varphi_c(x)\, dx.$$

We want to apply the dominated convergence theorem to the integral on the right-hand side of (5.4), with $h$ and $c$ kept fixed, $l \to \infty$ and *uniformly* in $v \in [-1, 1]$. □

*Choice of $\chi$.* At this point, we want to understand when (5.3) has a solution and we treat the cases $c \leq 1$ and $c > 1$ separately. If $c \leq 1$, it is easy to check that equation (5.3) has a solution for all $v \in [-1, 1]$ if and only if $h > c$ and moreover the solution is then unique for all such $v$.

The $c > 1$ case is a bit longer. Let $\chi^*$ be

(5.5) $$\chi^* = \chi^*(c) := \ln(\sqrt{c} + \sqrt{c-1}\,),$$

which is equivalent to

(5.6) $$(\cosh \chi^*)^2 = c \quad \text{and} \quad \chi^* > 0,$$

or to the fact that $\chi^*$ is the unique local minimum of $\xi \mapsto \xi - c \tanh \xi$. Observe that

(5.7) $$\beta(c) := \ln(\sqrt{c} + \sqrt{c-1}\,) + c - \sqrt{c^2 - c} = \chi^* - c \tanh \chi^* + c.$$

This is, of course the same as

(5.8) $$\beta(c) - c = \min_{\xi \geq 0}(\xi - c \tanh \xi).$$

Again, (5.3) has a solution for all $v \in [-1, 1]$ if and only if

(5.9) $$h \geq \beta(c).$$

For $h > \beta(c)$, we choose the solution of (5.3) with

$$\chi > \chi^*.$$

Given $c$ and $h > \beta(c)$, we denote by $\bar\chi = \bar\chi(c, h)$ the solution of the equation (5.3) with $v = +1$. Clearly,

$$\chi^* < \bar\chi = \min_{-1 \leq v \leq +1} \chi$$

and hence, due to (5.6) we always have

(5.10) $$\sup_{-1 \leq v \leq +1} 4pqc = 4\bar p \bar q c < 4p^* q^* c = 1,$$

where $\bar p, \bar q, p^*$ and $q^*$ all have the obvious meaning. Furthermore, keeping $c$ fixed, $[\beta(c), \infty) \ni h \mapsto \bar\chi(c, h) \in [\chi^*(c), \infty)$ is strictly increasing in $h$ with

(5.11) $$\lim_{h \to \infty} \bar\chi(c, h) = \infty.$$

*Pointwise convergence.* For fixed $c$ and $h > \beta(c)$, the integrand on the right-hand side of (5.4) converges pointwise to

(5.12) $$\frac{\exp\{-x^2(1 - 4pqc)/(2c)\}}{\sqrt{2\pi c}},$$

as $l \to \infty$ uniformly on compact domains of $x$ and in the parameter $v \in [-1, 1]$. Due to (5.10), which holds for all $c$ and $h > \beta(c)$, the limit function is integrable *uniformly* for $v \in [-1, 1]$.

*Domination.* For $\varepsilon \in [0,1]$ we define

$$\tilde{c}(\varepsilon) := \sup\left\{c : \inf_{-\infty < y < \infty} e^{y^2/c}(1 - \varepsilon(\sin y)^2) = 1\right\}. \tag{5.13}$$

Note that $\varepsilon \mapsto \tilde{c}(\varepsilon)$ is monotone decreasing, with

$$\lim_{\varepsilon \searrow 0} \tilde{c}(\varepsilon) = \infty \quad \text{and} \quad \tilde{c}(1) := \lim_{\varepsilon \nearrow 1} \tilde{c}(\varepsilon) = 0. \tag{5.14}$$

It follows that for all $c, h$ and $v$

$$\inf_{-\infty < y < \infty} |e^{y^2/(2\tilde{c}(4pq))}(pe^{i2qy} + qe^{-i2py})| = 1, \tag{5.15}$$

and this will be used in order to bound the integrand on the right-hand side of (5.4).

LEMMA 5.1. *We have:*

(i) *For any $\varepsilon \in [0,1]$*

$$-\frac{1}{\ln(1-\varepsilon)} \leq \tilde{c}(\varepsilon) \leq \min\left\{-\frac{\pi^2}{4\ln(1-\varepsilon)}, \frac{1}{\varepsilon}\right\}. \tag{5.16}$$

(ii) *For $\varepsilon \leq 2/3$*

$$\tilde{c}(\varepsilon) = \frac{1}{\varepsilon}. \tag{5.17}$$

PROOF. (i) We obtain the first upper bound in (5.16) by looking at $y = \pi/2$ and we obtain the second upper bound by expanding near $y = 0$ in (5.13). In order to prove the lower bound of (5.16) note that for any $\varepsilon, \alpha \in [0,1]$,

$$1 - \varepsilon\alpha \geq (1-\varepsilon)^\alpha.$$

Hence

$$1 - \varepsilon(\sin y)^2 \geq \exp\{(\sin y)^2 \ln(1-\varepsilon)\} \geq \exp\{y^2 \ln(1-\varepsilon)\}.$$

(ii) In view of the first part of this lemma, we need only check that for $\varepsilon \leq 2/3$,

$$g(y) = e^{\varepsilon y^2}(1 - \varepsilon \sin^2 y) \geq 1, \qquad y \in \mathbb{R}.$$

To do so, note that $g(0) = 1$ and compute

$$g'(y) = 2\varepsilon e^{\varepsilon y^2} h(y),$$

where

$$h(y) = y(1 - \varepsilon \sin^2(y)) - \sin y \cos y.$$

Then
$$h'(y) = \sin y((2-\varepsilon)\sin y - 2\varepsilon y \cos y),$$

which is zero if (i) $\sin y = 0$ or (ii)
$$y = \frac{(2-\varepsilon)\sin y}{2\varepsilon \cos y}.$$

In case (i), $h(y) = y$, while in case (ii),
$$h(y) = \frac{\sin y(2 - 3\varepsilon + \varepsilon^2 \sin^2 y)}{2\varepsilon \cos y} = \frac{y(2 - 3\varepsilon + \varepsilon^2 \sin^2 y)}{2 - \varepsilon}.$$

If $\varepsilon \leq 2/3$, we see that $h(y)$ and $y$ have the same sign at each critical point of $h$. Since $h(y) \to \infty$ as $y \to \infty$ in this case, it follows that $h(y) \geq 0$ for $y \geq 0$, and hence that $g$ is increasing on $[0,\infty)$. Therefore, $g(y) \geq 1$ for all $y$. □

Now we return to the boundedness of the integrand on the right-hand side of (5.4). Let

(5.18)
$$\bar{c} = \bar{c}(c,h) := \min_{-1 \leq v \leq +1} \tilde{c}(4pq)$$
$$= \tilde{c}\left(\max_{-1 \leq v \leq +1} 4pq\right) = \tilde{c}(1/(\cosh \bar{\chi}(c,h))^2).$$

LEMMA 5.2. (i) If

(5.19) $$c < \bar{c}(c,h)$$

holds, then the integrand on the right-hand side of (5.4) is bounded by $e^{-x^2(\bar{c}-c)/(2\bar{c}c)}/\sqrt{2\pi c}$ for all $x, l$ and $v \in [-1,1]$.

(ii) For any $c$, if $h > h^*(c)$, then (5.19) holds.

PROOF. (i) We clearly have
$$|(pe^{ix/\sqrt{l}} + qe^{-ix/\sqrt{l}})^n| \leq 1.$$

Using $\bar{c}$ defined in (5.18) we write
$$(pe^{i2qx/\sqrt{l}} + qe^{-i2px/\sqrt{l}})^{-l} e^{-x^2/(2c)}$$
$$= (e^{(x/\sqrt{l})^2/(2\bar{c})}(pe^{i2qx/\sqrt{l}} + qe^{-i2px/\sqrt{l}}))^{-l} e^{-x^2(\bar{c}-c)/(2c\bar{c})}.$$

From (5.15) and (5.18) it follows that the absolute value of this last expression is bounded by $e^{-x^2(\bar{c}-c)/(2c\bar{c})}$.

(ii) Due to (5.11) and (5.14), it is clear that for any $c$, (5.19) holds for all large $h$. Carrying out a tedious calculation leads to the statement for

$c < 3/2$. For the case $c \geq 3/2$, we need only observe that once $h > \beta(c)$, we have that

$$2/3 \geq \frac{1}{c} = \frac{1}{(\cosh \chi^*(c))^2} > \frac{1}{(\cosh \bar{\chi}(c,h))^2} = \frac{1}{\tilde{c}(1/(\cosh \bar{\chi}(c,h))^2)},$$

the last equality following from (5.17). $\square$

CONCLUSION OF THE PROOF OF THEOREM 1.21. Fix $c > 0$ and $h > h^*(c)$. Lemma 5.2(ii) tells us that (5.19) holds. It then follows from the uniform convergence in (5.12) and Lemma 5.2(i) that the integral on the right-hand side of (5.4) converges to $1/\sqrt{1-4pqc}$, as $l \to \infty$, uniformly in $v \in [-1,1]$. Since $1/\sqrt{1-4pqc}$ is clearly bounded away from 0 uniformly in $v$, it follows that this integral is positive for all large $l$. $\square$

**Acknowledgments.** Having been unaware of [20], in the first version of this paper, we had two proofs of Theorem 1.1, the first one turning out to follow the lines of that in [20] and the other being the present alternative proof. We thank Yuval Peres for bringing [20] to our attention. In addition, in the first version of this paper, we had Proposition 4.2 but were unaware of it having been known. We thank Sarbarish Chakravarty for bringing [9] to our attention. We also thank David Brydges for bringing certain references to our attention.

T. M. LIGGETT
DEPARTMENT OF MATHEMATICS
UNIVERSITY OF CALIFORNIA
   AT LOS ANGELES
LOS ANGELES, CALIFORNIA 90095
USA
E-MAIL: tml@math.ucla.edu

J. E. STEIF
DEPARTMENT OF MATHEMATICS
CHALMERS UNIVERSITY OF TECHNOLOGY
   AND GÖTEBORG UNIVERSITY
S-412 96 GOTHENBURG
SWEDEN
E-MAIL: steif@math.chalmers.se

B. TÓTH
INSTITUTE OF MATHEMATICS
TECHNICAL UNIVERSITY BUDAPEST
EGRY JÓZSEF U. 1
II-111 BUDAPEST
HUNGARY
E-MAIL: balint@math.bme.hu